# PERCOLATION ON DENSE GRAPH SEQUENCES

By Béla Bollobás[1], Christian Borgs, Jennifer Chayes
and Oliver Riordan[2]

*University of Cambridge and University of Memphis, Microsoft Research,
Microsoft Research and University of Oxford*

In this paper we determine the percolation threshold for an arbitrary sequence of dense graphs $(G_n)$. Let $\lambda_n$ be the largest eigenvalue of the adjacency matrix of $G_n$, and let $G_n(p_n)$ be the random subgraph of $G_n$ obtained by keeping each edge independently with probability $p_n$. We show that the appearance of a giant component in $G_n(p_n)$ has a sharp threshold at $p_n = 1/\lambda_n$. In fact, we prove much more: if $(G_n)$ converges to an irreducible limit, then the density of the largest component of $G_n(c/n)$ tends to the survival probability of a multi-type branching process defined in terms of this limit. Here the notions of convergence and limit are those of Borgs, Chayes, Lovász, Sós and Vesztergombi.

In addition to using basic properties of convergence, we make heavy use of the methods of Bollobás, Janson and Riordan, who used multi-type branching processes to study the emergence of a giant component in a very broad family of sparse inhomogeneous random graphs.

**1. Introduction.** In this paper we study percolation on arbitrary sequences of dense finite graphs, where the number of edges grows quadratically with the number of vertices. The study of percolation on finite graphs is much more delicate than that of percolation on infinite graphs; indeed, percolation on finite graphs provides the finite-size scaling behavior of percolation on the corresponding infinite graphs (see, e.g., Borgs, Chayes, Kesten and Spencer [11] for the study of percolation on finite subcubes of $\mathbb{Z}^d$).

Received February 2007; revised April 2008.
[1]Supported in part by NSF Grants DMS-0505550, CNS-0721983 and CCF-0728928 and ARO Grant W911NF-06-1-0076.
[2]Supported by a Royal Society research fellowship held at the Department of Pure Mathematics and Mathematical Statistics, Cambridge CB3 0WB, UK.

*AMS 2000 subject classifications.* 60K35, 05C80.
*Key words and phrases.* Percolation, cut metric, random graphs.

The first question one asks is whether there is a percolation phase transition. In the case of a finite graph on $n$ vertices, we say that a percolation phase transition occurs when the size of the largest component goes from being of order $o(n)$ [typically in fact $O(\log n)$] below a certain density to order $n$ above that density. Of course, to make this precise, one must consider a sequence of graphs with $n \to \infty$. The next question one typically asks is how the size of the second largest component behaves. In the few specific cases studied so far, the second largest component is of order $\log n$ both below and above the transition; the behavior above the transition is much more difficult to prove. Once the existence of the transition has been established, one then studies the finite-size scaling (i.e., behavior in $n$) of the width of the transition region, and the size of the largest component within that transition window.

In this paper we establish the existence of a phase transition, including the behavior of the largest and second largest components, for a very large class of sequences of dense finite graphs. Moreover, we establish the location of the transition in terms of spectral properties of these graphs.

Consider a sequence of dense graphs $(G_n)$ and a sequence of random subgraphs $G_n(p_n)$ obtained from $G_n$ by deleting edges independently with probability $1 - p_n$. We say that the system *percolates* if $G_n(p_n)$ has a giant component, that is, a connected component of size $\Theta(|G_n|)$, where $|G_n|$ denotes the number of vertices in $G_n$. As usual, we say that the appearance of a giant component has a *sharp threshold* if there exists a sequence $(p_n)$ such that for all $\varepsilon > 0$, the random subgraph $G_n(p_n(1-\varepsilon))$ has no giant component with probability $1 - o(1)$ while $G_n(p_n(1+\varepsilon))$ has a giant component with probability $1 - o(1)$. (Here, and throughout, all asymptotic notation refers to the limit as $n \to \infty$.)

The simplest sequence $(G_n)$ for which this question has been analyzed is a sequence of complete graphs on $n$ vertices. The corresponding random subgraph is the well-known random graph $G_{n,p_n}$. Erdős and Rényi [16] were the first to show that with $p_n = c/n$, the random graph $G_{n,c/n}$ undergoes a phase transition at $c = 1$: for $c < 1$, all components are of size $O(\log n)$ while for $c > 1$ a giant component of size $\Theta(n)$ emerges. Later, the precise window of this phase transition was determined by Bollobás [5] and Łuczak [19].

Other specific sequences were considered in both the combinatorics and the probability communities. Ajtai, Komlós and Szemerédi [1] established a phase transition for percolation on the $n$-cube $Q_n = \{0,1\}^n$ (see [7, 10] for much more detailed estimates on this transition). Borgs, Chayes, Kesten and Spencer [11] studied the case when the graphs $G_n$ are rectangular subsets of $\mathbb{Z}^2$, and determined both the width of the phase transition window and the size of the largest component within this window in terms of the critical exponents of the infinite graph $\mathbb{Z}^2$.

While the question of a phase transition for random subgraphs of general sequences $(G_n)$ was already formulated by Bollobás, Kohayakawa and Łuczak [7], progress on this question has been rather slow. The few papers which deal with more general classes of graph sequences are still restricted in scope. See, for example, Borgs, Chayes, van der Hofstad, Slade and Spencer [8, 9] where the window for transitive graphs obeying the so-called triangle condition was analyzed, Frieze, Krivelevich and Martin [17] where the threshold for random subgraphs of a sequence of quasi-random graphs was analyzed, and Alon, Benjamini and Stacey [2] for results about expander graphs with bounded degrees.

Here we analyze the phase transition for random subgraphs of dense convergent graph sequences. The concept of convergent graph sequences was introduced for sparse graphs by Benjamini and Schramm [3] and for dense graphs by Borgs, Chayes, Lovász, Sós and Vesztergombi in [12] (see also [13]). As shown in [14, 15], there are many natural, a priori distinct definitions of convergence which turn out to be equivalent. Here we use the following one: given two graphs $F$ and $G$, define the *homomorphism density*, $t(F,G)$, of $F$ in $G$ as the probability that a random map from the vertex set of $F$ into the vertex set of $G$ is a homomorphism; a sequence $(G_n)$ of graphs is then said to be *convergent* if $t(F,G_n)$ converges for all finite graphs $F$. Note that any sequence $(G_n)$ has a convergent subsequence, so when studying general sequences of graphs, we may as well assume convergence.

It was shown by Lovász and Szegedy [18] that if a graph sequence converges, then the limiting homomorphism densities can be expressed in terms of a measurable function, $W:[0,1]^2 \to [0,1]$, which can therefore be thought of as the limit of the graph sequence. Following [14], we call such functions and their generalizations *graphons*. More precisely, a graphon is a bounded measurable function $W:[0,1]^2 \to [0,\infty)$ with $W(x,y) = W(y,x)$ for all $x,y$. In [15] it is also shown that convergence implies convergence of the normalized spectra of the adjacency matrices to the spectrum of the limiting graphon considered as an operator on $L^2([0,1])$.

Independently of the results above, Bollobás, Janson and Riordan [6] introduced a very general model of inhomogeneous random graphs with bounded average degree, defined in terms of so-called *kernels*. Although kernels are reminiscent of graphons, they are more general; in particular, they can be unbounded. One of the aims of [6] was to prove precise results about the emergence of the giant component in a general class of random graphs.

Our main result in this paper says that a convergent graph sequence has a sharp percolation threshold, and, moreover, if the limiting graphon $W$ is irreducible, then the density of the largest component is asymptotically equal to the survival probability of a certain multi-type branching process defined in terms of $W$ (see Theorem 1 below).

As a corollary of this result, we obtain that the appearance of a giant component in an arbitrary sequence of dense graphs $(G_n)$ (convergent or not) has a sharp threshold at $p_n = 1/\lambda_n$, where $\lambda_n$ is the largest eigenvalue of the adjacency matrix of $G_n$ (see Theorem 2 below). As usual, a sequence $(G_n)$ with $|G_n| \to \infty$ is called *dense* if the average degree in $G_n$ is of order $\Theta(|G_n|)$.

To state our results precisely, we need some notation. Given two graphs $F$ and $G$, write $\hom(F, G)$ for the number of homomorphisms (edge-preserving maps) from $F$ into $G$, and

$$t(F, G) = |G|^{-|F|} \hom(F, G)$$

for its normalized form, the *homomorphism density*. Following [13], we call a sequence $(G_n)$ of graphs *convergent* if $t(F, G_n)$ converges for every graph $F$. It was shown in [18] that a sequence $(G_n)$ is convergent if and only if there exists a symmetric, Lebesgue-measurable function $W : [0, 1]^2 \to [0, 1]$ such that

(1) $\qquad t(F, G_n) \to t(F, W) \qquad$ for every graph $F$,

where

$$t(F, W) = \int_{[0,1]^{V(F)}} \prod_{ij \in E(F)} W(x_i, x_j) \prod_{i \in V(F)} dx_i$$

is the *homomorphism density* of $F$ in $W$. In this case the sequence is said to *converge* to $W$; in notation, $G_n \to W$.

We also need the notion of a weighted graph. For the purposes of this paper, a *weighted graph* $G$ on a vertex set $V$ is a symmetric function $\beta : (v, w) \mapsto \beta_{vw}$ from $V \times V$ to $[0, \infty)$ with $\beta_{vv} = 0$ for every $v \in V$. (We thus do not allow vertex weights, and also restrict ourselves to nonnegative edge weights, instead of the more general case of real-valued edge weights considered in [14, 15].) Graphs correspond naturally to weighted graphs taking values in $\{0, 1\}$. The definitions of $t(F, G)$ and of convergence extend naturally to weighted graphs: if $F$ is a graph on $[k]$ and $G$ is a weighted graph on $V$, then

$$t(F, G) = |G|^{-|F|} \sum_{v_1, \ldots, v_k \in V} \prod_{ij \in E(F)} \beta_{v_i v_j}.$$

Let $(G_n)$ be a sequence of weighted graphs. We write $\beta_{ij}(n)$ for the weight of the edge $ij$ in $G_n$, suppressing the dependence on $n$ when this does not lead to confusion. We shall assume throughout that $|G_n| \to \infty$, that $\beta_{\max} = \sup_{i,j,n} \beta_{ij}(n) < \infty$ and that $(G_n)$ is convergent in the sense that $t(F, G_n)$ converges for all unweighted graphs $F$, although we shall remind the reader of this assumption in key places. As shown in [14], the results of [18] immediately generalize to such sequences, implying the existence of a graphon $W$ such that $G_n \to W$ in the sense of (1). Let us note that in the context of unweighted graphs, graphons are usually defined to take values in $[0,1]$. For weighted graphs with unbounded edge-weights, one could also consider unbounded limit functions $W$; we do not consider this extension here. One can also consider signed edge weights, in which case the appropriate limit objects are *signed graphons*, that is, (bounded) symmetric measurable functions from $[0,1]^2$ to $\mathbb{R}$.

Given a signed graphon $W$, the *cut norm* of $W$ is

$$\|W\|_\square = \sup_{A,B \subset [0,1]} \left| \int_{A \times B} W(x,y)\,dx\,dy \right|,$$

where the supremum is taken over all measurable sets. A *rearrangement* $W^\phi$ of a graphon $W$ is the graphon defined by $W^\phi(x,y) = W(\phi(x), \phi(y))$ where $\phi : [0,1] \to [0,1]$ is a measure-preserving bijection. Finally, the *cut metric* is the pseudometric on graphons defined by

$$\delta_\square(W_1, W_2) = \inf_\phi \|W_1 - W_2^\phi\|_\square.$$

It is proved in [14] that $G_n \to W$ if and only if

$$\delta_\square(W_{G_n}, W) \to 0, \tag{2}$$

where $W_{G_n}$ denotes the piecewise-constant graphon naturally associated to the weighted graph $G_n$.

Given a weighted graph $G$, let $G(p)$ be the random graph on $V(G)$ in which edges are present independently, and the probability that $ij$ is an edge is $\min\{p\beta_{ij}, 1\}$. We shall lose nothing by assuming that $p\beta_{ij} < 1$, so we shall often write $p\beta_{ij}$ for $\min\{p\beta_{ij}, 1\}$. Alternatively, as in [6], we could take $p\beta_{ij}$ to be a Poisson "edge intensity," so the probability of an edge is $1 - \exp(-p\beta_{ij})$; this makes no difference to our results. Our aim is to study the giant component in the sequence $G_n(c/|G_n|)$. To do this, we shall consider a certain branching process associated to the graphon $cW$.

As in [6], there is a natural way to associate a multi-type branching process $\mathfrak{X}_W$ to a measurable $W : [0,1]^2 \to \mathbb{R}^+$: each generation consists of a finite set of particles with "types" in $[0,1]$. Given generation $t$, each particle in generation $t$ has children in the next generation independently of the other particles and of the history. If a particle has type $x$, then the types of its children are distributed as a Poisson process on $[0,1]$ with intensity measure $W(x,y)\,dy$, where $dy$ denotes the Lebesgue measure. In other words, the number of children with types in a measurable set $A \subset [0,1]$ is Poisson with mean $\int_{y \in A} W(x,y)\,dy$, and these numbers are independent for disjoint sets $A$. The first generation of $\mathfrak{X}_W$ consists of a single particle whose type $x$ is uniformly distributed on $[0,1]$. Often we consider the same branching process but started with a particle of a fixed type $x$: we write $\mathfrak{X}_W(x)$ for this process.

If, as we shall always assume in this paper, $W$ is bounded, then

$$\lambda(x) = \int_0^1 W(x,y)\,dy,$$

the expected number of children of a particle of type $x$, is bounded by $\|W\|_\infty$; in particular, this expected number is finite. Thus every particle always has a finite number of children, and the total size of $\mathfrak{X}_W$ is infinite if and only if the process $\mathfrak{X}_W$ survives for ever.

Writing $|\mathfrak{X}_W|$ for the total number of particles in all generations of the branching process $\mathfrak{X}_W$, let $\rho(W) = \mathbb{P}(|\mathfrak{X}_W| = \infty)$ be the "survival probability" of $\mathfrak{X}_W$, and let $\rho(W;x)$ be the survival probability of $\mathfrak{X}_W(x)$, the process started with a particle of type $x$. From basic properties of Poisson processes we have

$$\rho(W;x) = 1 - \exp\left(-\int W(x,y)\rho(W;y)\,dy\right) \tag{3}$$

for every $x$, and from the definitions of $\mathfrak{X}_W(x)$ and $\mathfrak{X}_W$ we have

$$\rho(W) = \int \rho(W;x)\,dx.$$

In general, the functional equation (3) may have several solutions. It is a standard result of the theory of branching processes (proved in the present context in [6], for example) that $\rho(W;x)$ is given by the largest solution, that is, the pointwise supremum of all solutions.

Let $T_W$ be the integral operator defined by

$$(T_W(f))(x) = \int W(x,y)f(y)\,dy.$$

From Theorem 6.1 of [6], we have $\rho(W) > 0$ if and only if $\|T_W\| > 1$, where $\|T_W\|$ is the $L^2$-norm of $T_W$.

We shall show that the condition $G_n \to W$ is strong enough to ensure that the branching process captures enough information about the graph to determine the asymptotic size of the giant component in $G_n(c/|G_n|)$. For this we need one more definition, corresponding roughly to connectedness.

A graphon $W$ is *reducible* if there is a measurable $A \subset [0,1]$ with $0 < \mu(A) < 1$ such that $W(x,y) = 0$ for almost every $(x,y) \in A \times A^c$. Otherwise, $W$ is *irreducible*. Using the equivalent condition (2) for convergence, together with Szemerédi's lemma, it is easy to show that if $G_n \to W$ with $W$ reducible, then the vertex set of each $G_n$ may be partitioned into two classes so that the induced graphs $H_n$ and $K_n$ converge to appropriate graphons $W_1$ and $W_2$, with $o(|G_n|^2)$ edges of $G_n$ joining $H_n$ to $K_n$. In other words, the graphs $G_n$ may be written as "almost disjoint" unions of convergent sequences $H_n$ and $K_n$. In the light of this observation, it will always suffice to consider the case where $W$ is irreducible.

Henceforth, for notational convenience, we shall assume that $G_n$ has $n$ vertices; we shall often take the vertex set to be $[n] = \{1, 2, \ldots, n\}$. Note that we do not require $G_n$ to be defined for every $n$: all results will extend trivially to the general case $|G_n| \to \infty$ by considering subsequences. Let $(X_n)$ be a sequence of nonnegative random variables and $(a_n)$ a sequence of nonnegative reals. As usual, we say that a sequence of events $\mathcal{E}_n$ holds *with high probability*, or *whp*, if $\mathbb{P}(\mathcal{E}_n) \to 1$ as $n \to \infty$. As in [6], we write $X_n = o_{\mathrm{p}}(a_n)$ if $X_n/a_n$ converges to 0 in probability, $X_n = O(a_n)$ whp if there is a constant $C$ such that $X_n \leq Ca_n$ holds whp, and $X_n = \Theta(a_n)$ whp if there are constants $0 < C_1 \leq C_2$ such that $C_1 a_n \leq X_n \leq C_2 a_n$ holds whp. Note that the assertion $X_n = O(a_n)$ whp is stronger than the assertion that $X_n/a_n$ is bounded in probability [sometimes written $X_n = O_{\mathrm{p}}(a_n)$].

THEOREM 1. *Let $(G_n) = (\beta_{ij}(n))_{i,j \in [n]}$ be a sequence of weighted graphs with $|G_n| = n$ and $\beta_{\max} = \sup_{i,j,n} \beta_{ij}(n) < \infty$ converging to a graphon $W$. Let $c > 0$ be a constant, and let $C_1 = C_1(n)$ denote the number of vertices in the largest component of the random graph $G_n(c/n)$, and $C_2 = C_2(n)$ the number of vertices in the second largest component.*

*(a) If $c \leq \|T_W\|^{-1}$, then $C_1 = o_{\mathrm{p}}(n)$.*

*(b) If $c > \|T_W\|^{-1}$, then $C_1 = \Theta(n)$ whp. More precisely, for any constant $\alpha < (c\|T_W\| - 1)/(c\beta_{\max})$ we have $C_1 \geq \alpha n$ whp.*

*(c) If $W$ is irreducible, then $C_1/n \xrightarrow{\mathrm{p}} \rho(cW)$ and $C_2 = o_{\mathrm{p}}(n)$.*

In the result above, we may take $c = 1$ without loss of generality, rescaling the edge weights in $G_n$. The heart of the theorem is part (c); as we shall see later [in the discussion around (14)], part (b) follows easily.

The first two statements of Theorem 1 immediately imply that an arbitrary sequence of dense graphs has a sharp percolation threshold.

THEOREM 2. *Let $(G_n)$ be a sequence of dense graphs with $|G_n| = n$, let $\lambda_n$ be the largest eigenvalue of the adjacency matrix of $G_n$ and let $p_n = \min\{c/\lambda_n, 1\}$.*

*(a) If $c \leq 1$, then the largest component of $G_n(p_n)$ is of size $o_{\mathrm{p}}(n)$.*

*(b) If $c > 1$, then the largest component of $G_n(p_n)$ has size $\Theta(n)$ whp.*

PROOF. Theorem 2 follows immediately from Theorem 1 and the fact that any sequence of weighted graphs with uniformly bounded edge weights has a convergent subsequence. Indeed, to see (a), suppose for a contradiction that $c \leq 1$, and for some $\varepsilon > 0$ there is a subsequence (which we again denote by $G_n$) such that with probability at least $\varepsilon$, the largest component of $G_n(p_n)$ has at least $\varepsilon n$ vertices. Let $\tilde{G}_n$ be the sequence of weighted graphs obtained by weighting each edge in $G_n$ by $n/\lambda_n$. Since $(G_n)$ is dense, we have $\lambda_n = \Theta(n)$, so the edge weights in $(\tilde{G}_n)$ are bounded above. By compactness, we may assume that the homomorphism densities $t(F, \tilde{G}_n)$ converge, so by the results of [14], we may assume that $\tilde{G}_n$ converges [in the sense of (1)] to some graphon $W$. Since the largest eigenvalue of the adjacency matrix of $\tilde{G}_n$ is equal to $n$, we have that $\|T_W\| = 1$, so by Theorem 1 the largest component of $G_n(p_n) = \tilde{G}_n(c/n)$ has size $o_p(n)$, a contradiction. The proof of (b) proceeds along the same lines. □

For convergent sequences of weighted graphs converging to an irreducible graphon $W$, we shall prove stronger results about the sizes of the small components in the noncritical cases. This time we renormalize by scaling the edge weights, taking $c = 1$.

THEOREM 3. *Let $(G_n)$ be a sequence of edge-weighted graphs with $|G_n| = n$ and $\beta_{\max} = \sup_{i,j,n} \beta_{ij}(n) < \infty$, and suppose that $G_n \to W$.*

*If $\|T_W\| < 1$, then there is a constant $A$ such that $C_1(G_n(1/n)) \le A \log n$ holds whp.*

*If $\|T_W\| > 1$ and $W$ is irreducible, then there is a constant $A$ such that $C_2(G_n(1/n)) \le A \log n$ holds whp.*

It is easily seen that $\|T_W\| > 1$ does not in itself imply that the second component has $O(\log n)$ vertices; indeed, for an appropriate sequence $(G_n)$ in which each $G_n$ is the disjoint union of two graphs, $H_n$ and $K_n$, with $(H_n)$ and $(K_n)$ convergent, $H_n(1/n)$ and $K_n(1/n)$ both have giant components whp, so the second component of $G_n(1/n)$ has $\Theta(n)$ vertices whp. Note also that adding, say, $\Theta(n^{3/2})$ random edges to $G_n$ running from $H_n$ to $K_n$ will almost certainly join [in $G_n(1/n)$] any giant components in $H_n(1/n)$ and $K_n(1/n)$ while preserving the condition $G_n \to W$. Hence we cannot expect to find the asymptotic size of the giant component in the reducible case.

Let $N_k(G)$ denote the number of vertices of a graph $G$ that are in components of size (number of vertices) exactly $k$. The basic idea of the proof of Theorem 1 is to consider components of each fixed size; the key lemma we shall need is as follows.

LEMMA 4. *Let $W:[0,1]^2 \to [0,\infty)$ be a graphon, and let $(G_n) = (\beta_{ij}(n))_{i,j \in [n]}$ be a sequence of weighted graphs with $G_n \to W$ and $\beta_{\max} = \sup_{i,j,n} \beta_{ij}(n)$ finite. For each fixed $k$ we have*

$$\frac{1}{n} N_k(G_n(1/n)) \overset{\mathrm{p}}{\to} \mathbb{P}(|\mathfrak{X}_W| = k) \tag{4}$$

*as $n \to \infty$.*

Lemma 4 tells us that we have the "right" number of vertices in small components; we shall then show that most of the remaining vertices are in a single large component. Of course, as with any such branching process lemma, the proof actually gives a little more: for any finite tree $T$, $1/n$ times the number of tree components of $G_n(1/n)$ isomorphic to $T$ converges in probability to a quantity that may be calculated from the branching process. In particular, considering a rooted tree $T^*$, the normalized number of vertices $v$ in tree components that are isomorphic to $T^*$ with $v$ as the root converges to the probability that the branching process $\mathfrak{X}_W$ is isomorphic to $T^*$, when $\mathfrak{X}_W$ is viewed as a rooted tree in the natural way.

Lemma 4 will be proved in the next section; the proof of Theorem 1 is then given in Section 3. Theorem 3 is proved in Section 4 in two separate parts.

Before turning to the proofs, let us remark on the relationship of the present results to those of [6]. It might at first sight appear that Theorem 1 (and hence Theorem 2) would follow simply from Theorem 3.1 of [6]. Indeed, given a sequence $(G_n)$, the sequence of random subgraphs $G_n(c/n)$ can be seen as an instance of the model of [6] with a sequence $\kappa_n$ of kernels given by $\kappa_n = cW_{G_n}$. Passing to a subsequence, we can assume that these kernels converge in the cut metric. But the results of [6] require pointwise convergence almost everywhere, which is much stronger. Considering in particular the case where the $G_n$ are unweighted graphs, the kernels $\kappa_n$ take only the values $c$ and $0$, so it is clear that in general they will have no subsequence that converges pointwise.

One could of course ask whether Theorem 3.1 of [6] may be strengthened by relaxing the pointwise convergence condition to convergence in the cut metric. Ignoring complications due to the more general setting considered in [6], the answer for bounded kernels turns out to be yes, but that is essentially the content of the present paper. Indeed, as noted above, Theorem 1 may be viewed as a very special case of such a strengthening, for kernels of a special form (piecewise constant on each $1/n$ by $1/n$ square) and the corresponding vertex space. It is easy to extend this to general sequences of uniformly bounded kernels as long as, in the terminology of [6], we restrict ourselves to the vertex space in which the types of the vertices are independent and uniformly distributed on $[0,1]$. We can state such an extension without reference to the more complicated definitions in [6].

Given a graphon $W$, for $n \geq \sup W$ let $G(n, W)$ be the random (sparse) graph on $[n]$ obtained as follows: first choose $X_1, \ldots, X_n$ independently and uniformly from $[0,1]$. Then, given $X_1, \ldots, X_n$, join each pair $ij$, $i < j$, with probability $W(X_i, X_j)/n$, independently of the other pairs.

THEOREM 5. *Let $W_n$ be a sequence of uniformly bounded graphons and $W$ a graphon, and suppose that $\delta_\square(W_n, W) \to 0$. Then the conclusions of Theorem 1 hold with $G_n(c/n)$ replaced by $G(n, cW_n)$.*

PROOF. This is essentially immediate from Theorem 4.7 of [14] and Theorem 1. Indeed, taking $X_1, \ldots, X_n$ independent and uniform on $[0,1]$, let $H_n = H(n, W_n)$ be the (dense) weighted graph on $[n]$ with weights $\beta_{ij}(n) = W_n(X_i, X_j)$. Theorem 4.7 of [14] states that with probability at least $1 - e^{-n^2/(2\log_2 n)}$ we have

$$\delta_\Box(W_n, W_{H_n}) \leq \frac{10}{\sqrt{\log_2 n}} \sup W_n \leq \frac{10M}{\sqrt{\log_2 n}},$$

where $M = \sup_n \sup W_n < \infty$. It follows that $\delta_\Box(W_n, W_{H_n}) \to 0$ a.s. Since $\delta_\Box(W_n, W) \to 0$, we thus have $H_n \to W$ a.s. To construct $G(n, W_n)$, we condition on the sequence $H_n$, and then realize $G(n, W_n)$ as $H_n(c/n)$. In this conditioning we may assume that $H_n \to W$, so the result follows by Theorem 1. □

Note that while Theorem 5 greatly strengthens Theorem 3.1 of [6] by relaxing the convergence condition to one that can always be applied to a subsequence, it is much weaker in other ways: we must restrict to uniformly bounded graphons/kernels, and there is no obvious way to handle general vertex spaces.

**2. The number of vertices in small components.** Let us first prove a slightly weaker form of the special case $k = 2$ of Lemma 4, calculating the expected number of isolated edges in $G_n(1/n)$. The only extra complications in the general case will be notational. Let

$$d_v = \sum_w \beta_{vw}$$

denote the "weighted degree" of a vertex $v$ in the weighted graph $G_n$. Note that for $v$ fixed, the quantity $d_v = d_v(n)$ is determined by $G_n$, so it is deterministic.

Let $v$ and $w$ be two vertices of $G_n$ chosen independently and uniformly at random, independently of which edges are selected to form $G_n(1/n)$. Given a random variable $X$ that depends on $G_n$, $v$ and $w$, but not on which edges are selected to form $G_n(1/n)$, we shall write $\mathbb{E}_{vw}$ for the expectation of $X$ over the choice of $v$ and $w$. We define $\mathbb{E}_v$ similarly.

Let $\mathbf{X}$ be the vector-valued random variable $\mathbf{X} = (d_v/n, d_w/n, \beta_{vw})$. Note that $\|\mathbf{X}\|_\infty \leq \sup_{v,w,n} \beta_{vw}(n) < \infty$. For nonnegative integers $t_1$ and $t_2$, let

(5) $$\mathbf{X}_{t_1, t_2} = \mathbb{E}_{vw}((d_v/n)^{t_1}(d_w/n)^{t_2} \beta_{vw})$$

denote the $(t_1, t_2, 1)$st joint moment of the triple $\mathbf{X}$. Then

$$n^{t_1 + t_2 + 2} \mathbf{X}_{t_1, t_2} = \sum_v \sum_w d_v^{t_1} d_w^{t_2} \beta_{vw}$$

$$= \sum_v \sum_w \left(\sum_u \beta_{vu}\right)^{t_1} \left(\sum_x \beta_{wx}\right)^{t_2} \beta_{vw}$$

$$= \sum_v \sum_w \sum_{u_1,\ldots,u_{t_1}} \sum_{x_1,\ldots,x_{t_2}} \beta_{vw} \prod_i \beta_{vu_i} \prod_i \beta_{wx_i}$$

$$= n^{t_1+t_2+2} t(S_{t_1,t_2}, G_n),$$

where $S_{t_1,t_2}$ is the "double star" consisting of an edge with $t_1$ extra pendent edges added to one end and $t_2$ to the other. Note that the summation variables are not required to be distinct, and that $t(F, G_n)$ counts homomorphisms, not injections. Since $G_n \to W$, it follows that

(6) $$\mathbf{X}_{t_1,t_2} \to t(S_{t_1,t_2}, W)$$

as $n \to \infty$.

For a given pair of vertices $v$, $w$, the probability [when we make the random choices determining $G_n(1/n)$] that $vw$ forms an isolated edge in $G_n(1/n)$ is exactly

$$p_{vw} = \frac{\beta_{vw}}{n} \prod_{z \neq v,w} \left(1 - \frac{\beta_{vz}}{n}\right)\left(1 - \frac{\beta_{wz}}{n}\right).$$

Note that the probability that $v$ is one end of an isolated edge is

$$\sum_{w \neq v} p_{vw} = n\mathbb{E}_w p_{vw},$$

so, with both $v$ and $w$ random, we have

$$\frac{1}{n}\mathbb{E}(N_2(G_n(1/n))) = \mathbb{P}(v \text{ in isol. edge}) = \mathbb{E}_{vw}(np_{vw}).$$

Since the $\beta$'s are bounded, we have

(7) $$p_{vw} \sim \frac{\beta_{vw}}{n} \exp\left(-\frac{d_v}{n} - \frac{d_w}{n}\right).$$

Let $D_{vw} = d_v/n + d_w/n$. Note that with $v$ and $w$ fixed, the quantity $D_{vw} = D_{vw}(n)$ is determined by $G_n$, so it is deterministic. For every $v$ and $w$ we have

$$\beta_{vw} \exp(-D_{vw}) = \sum_{t=0}^{\infty} (-1)^t \frac{\beta_{vw} D_{vw}^t}{t!}.$$

Taking $v$ and $w$ uniformly random, for each fixed $n$ we have

(8) $$\mathbb{E}_{vw}(\beta_{vw} \exp(-D_{vw})) = \sum_{t=0}^{\infty} (-1)^t \frac{\mathbb{E}_{vw}(\beta_{vw} D_{vw}^t)}{t!}.$$

As $n \to \infty$, from (5) and (6) we have

$$\mathbb{E}_{vw}(\beta_{vw} D_{vw}^t) = \mathbb{E}_{vw}(\beta_{vw}(d_v/n + d_w/n)^t)$$

$$= \sum_{t_1+t_2=t} \binom{t}{t_1} \mathbf{X}_{t_1,t_2} \to \sum_{t_1+t_2=t} \binom{t}{t_1} t(S_{t_1,t_2}, W).$$

It is easy to see that the $n \to \infty$ limit may be taken inside the sum in (8); indeed, $D_{vw}$ is bounded by $2\beta_{\max}$, so the $t$th summand in (8) is bounded by $C_t = \beta_{\max}(2\beta_{\max})^t/t!$. Since $\sum_t C_t = \beta_{\max} \exp(2\beta_{\max}) < \infty$, the sum in (8) is absolutely convergent, uniformly in $n$. Hence,

$$\lim_{n \to \infty} \mathbb{E}_{vw}(\beta_{vw} \exp(-D_{vw})) = \sum_{t=0}^{\infty} (-1)^t \frac{\lim_{n \to \infty} \mathbb{E}_{vw}(\beta_{vw} D_{vw}^t)}{t!}$$

$$= \sum_{t=0}^{\infty} \frac{(-1)^t}{t!} \sum_{t_1+t_2=t} \binom{t}{t_1} t(S_{t_1,t_2}, W).$$

Putting the pieces together, we have expressed the limiting expectation of $N_2(G_n(1/n))$ in terms of $W$:

$$\frac{1}{n}\mathbb{E}(N_2(G_n(1/n))) = \mathbb{E}_{vw}(np_{vw}) \sim \mathbb{E}_{vw}(\beta_{vw} \exp(-D_{vw}))$$

$$\to \sum_{t=0}^{\infty} \frac{(-1)^t}{t!} \sum_{t_1+t_2=t} \binom{t}{t_1} t(S_{t_1,t_2}, W)$$

$$= \sum_{t_1=0}^{\infty} \sum_{t_2=0}^{\infty} \frac{(-1)^{t_1}}{t_1!} \frac{(-1)^{t_2}}{t_2!} t(S_{t_1,t_2}, W).$$

Recalling that we write $\lambda(x) = \int_y W(x, y) \, dy$, we have

$$t(S_{t_1,t_2}, W) = \int_x \int_y W(x,y) \lambda(x)^{t_1} \lambda(y)^{t_2} \, dx \, dy,$$

so the final quantity above is simply

$$\sum_{t_1=0}^{\infty} \sum_{t_2=0}^{\infty} \int_x \int_y W(x,y) \frac{(-\lambda(x))^{t_1}}{t_1!} \frac{(-\lambda(y))^{t_2}}{t_2!} \, dx \, dy.$$

As $W$ is bounded, using dominated convergence we may take the sums inside the integral, obtaining

$$\int_x \int_y W(x,y) e^{-\lambda(x)} e^{-\lambda(y)} \, dx \, dy = \mathbb{P}(|\mathfrak{X}_W| = 2).$$

We have thus proved a weak form of (4) for $k = 2$, namely convergence in expectation: $\frac{1}{n}\mathbb{E}(N_2(G_n(1/n))) \to \mathbb{P}(|\mathfrak{X}_W| = 2)$.

Convergence in expectation for general $k$ is essentially the same, although we must count trees with each possible structure separately. Note that we need only consider tree components: if $N'_k(G)$ denotes the number of vertices of a graph $G$ that are in components of size $k$ that are not trees, then $N'_k(G)$ is certainly bounded by $k$ (or in fact $k/3$) times the number of vertices of $G$ in cycles of length at most $k$. In $G_n(1/n)$, this latter quantity has expectation at most

$$\sum_{\ell=3}^{k} n^\ell (\beta_{\max}/n)^\ell \le k \max\{1, \beta_{\max}^k\},$$

so we certainly have

(9) $$\mathbb{E}(N'_k(G_n(1/n))) \le k^2 \max\{1, \beta_{\max}^k\} = o(n)$$

as $n \to \infty$ with $k$ fixed. After this preparation, let us now turn to the proof of Lemma 4.

PROOF OF LEMMA 4. Let $T$ be a rooted tree on $k$ vertices. Let $\mathrm{aut}(T)$ denote the number of automorphisms of $T$ as a rooted tree. Thus, if $T_1, \ldots, T_r$ are the "branches" of $T$, then $\mathrm{aut}(T) = f \prod \mathrm{aut}(T_i)$, where the factor $f$ is the product of a factor $j!$ for each (maximal) set of $j$ isomorphic branches $T_i$.

The branching process $\mathfrak{X}_W$ may be naturally viewed as a rooted tree, by joining each particle to its parent and taking the initial particle as the root. We write $\mathfrak{X}_W \cong T$ if this tree is isomorphic to $T$ as a rooted tree. Note that

(10) $$\mathbb{P}(|\mathfrak{X}_W| = k) = \sum_T \mathbb{P}(\mathfrak{X}_W \cong T),$$

where the sum runs over all isomorphism classes of rooted trees on $k$ vertices.

Realizing $T$ as a graph on $\{1, 2, \ldots, k\}$, one can show that

(11) $$\mathbb{P}(\mathfrak{X}_W \cong T) = \frac{1}{\mathrm{aut}(T)} \int_{x_1} \cdots \int_{x_k} \prod_{i=1}^{k} e^{-\lambda(x_i)} \prod_{ij \in E(T)} W(x_i, x_j).$$

For example, the stronger statement,

$$\mathbb{P}(\mathfrak{X}_W(x_1) \cong T) = \frac{1}{\mathrm{aut}(T)} \int_{x_2} \cdots \int_{x_k} \prod_{i=1}^{k} e^{-\lambda(x_i)} \prod_{ij \in E(T)} W(x_i, x_j)$$

may be proved by induction on the size of $T$, noting that $\mathfrak{X}_W \cong T$ holds if and only if, for each isomorphism class of branch $T_i$ of $T$, we have exactly the right number of children of the initial particle whose descendants form a tree isomorphic to $T_i$.

Let $T$ be a tree on $\{1, 2, \ldots, k\}$, which we shall regard as a rooted tree with root 1, and let $\mathbf{v} = (v_1, \ldots, v_k)$ be a $k$-tuple of vertices of $G_n$. If $v_1, \ldots, v_k$ are distinct, let $p_{\mathbf{v},T}$ denote that probability that $\{v_1, \ldots, v_k\}$ is the vertex set of a component of $G_n(1/n)$ with $v_i v_j$ an edge of $G_n(1/n)$ if and only if $ij$ is an edge of $T$. (Note that this condition is stronger than the component being isomorphic to $T$.) If two or more $v_i$ are the same, set $p_{\mathbf{v},T} = 0$. As before, let $d_v = \sum_{w \in V(G_n)} \beta_{vw}$ denote the "degree" of $v$ in $G_n$; also, let

$$\beta_{\mathbf{v},T} = \prod_{ij \in E(T)} \beta_{v_i v_j}.$$

Arguing as for (7), since $\beta_{\max} < \infty$ we have

$$p_{\mathbf{v},T} \sim \frac{\beta_{\mathbf{v},T}}{n^{k-1}} \exp\left(-\sum_{i=1}^{k} \frac{d_{v_i}}{n}\right).$$

Furthermore, taking $v_1, \ldots, v_k$ uniformly random, and recalling from (9) that we may ignore nontree components, we have

$$\frac{1}{n}\mathbb{E}(N_k(G_n(1/n)))$$
$$= \sum_{T} \mathbb{P}(\text{random vertex } v \text{ is root of cpt isom. to } T) + o(1)$$
$$= \frac{1}{n} \sum_{T} \frac{1}{\operatorname{aut}(T)} \sum_{v_1, \ldots, v_k} p_{\mathbf{v},T} + o(1),$$

where the sum runs over all isomorphism classes of $k$-vertex rooted trees $T$, and the factor $1/\operatorname{aut}(T)$ accounts for the number of labelings of a tree component isomorphic (as a rooted tree) to $T$ with a given vertex $v_1$ as the root, which gives the number of distinct $k$-tuples $(v_1, \ldots, v_k)$ corresponding to a certain rooted component. Putting the above together, we have

$$\frac{1}{n}\mathbb{E}(N_k(G_n(1/n))) = n^{-k} \sum_{T} \sum_{\mathbf{v}} \frac{\beta_{\mathbf{v},T}}{\operatorname{aut}(T)} \exp\left(-\sum_{i=1}^{k} \frac{d_{v_i}}{n}\right) + o(1)$$
$$= \sum_{T} \mathbb{E}_{\mathbf{v}}\left(\frac{\beta_{\mathbf{v},T}}{\operatorname{aut}(T)} \exp\left(-\sum_{i=1}^{k} \frac{d_{v_i}}{n}\right)\right) + o(1),$$

where the expectation is over the uniformly random choice of $v_1, \ldots, v_k$.

The rest of the calculations are as before: it is easy to check that for nonnegative integers $t_1, \ldots, t_k$ we have

$$n^{\sum t_i + k} \mathbb{E}_{\mathbf{v}}\left(\beta_{\mathbf{v},T} \prod_{i=1}^{k} \left(\frac{d_{v_i}}{n}\right)^{t_i}\right) = n^{\sum t_i + k} t(T_{\mathbf{t}}, G_n),$$

where $T_{\mathbf{t}}$ is the graph formed from $T$ by adding $t_i$ pendent edges to each vertex $i$, so

$$t(T_{\mathbf{t}}, G_n) \to t(T_{\mathbf{t}}, W) = \int_{x_1} \cdots \int_{x_k} \prod_{ij \in E(T)} W(x_i, x_j) \prod_{i=1}^{k} \lambda(x_i)^{t_i}.$$

As all relevant sums are uniformly absolutely convergent, we can expand the term $\exp(-\sum d_{v_i}/n)$ as a sum of terms of the form $\prod_i (d_{v_i}/n)^{t_i}$, sum and take limits as before, finally obtaining

$$\frac{1}{n} \mathbb{E}(N_k(G_n(1/n))) = \sum_T \frac{1}{\operatorname{aut}(T)} \int_{x_1} \cdots \int_{x_k} \prod_{ij \in E(T)} W(x_i, x_j) \prod_{i=1}^{k} e^{-\lambda(x_i)} + o(1).$$

But, from (11), the final summand is just $\mathbb{P}(\mathfrak{X}_W \cong T)$, so, from (10),

$$\frac{1}{n} \mathbb{E}(N_k(G_n(1/n))) \to \mathbb{P}(|\mathfrak{X}_W| = k).$$

To complete the proof of Lemma 4 it remains to give a suitable upper bound on the variance. Let $N_{\geq k}(G)$ denote the number of vertices of a graph $G$ in components of size at least $k$, and set $N_{\geq k} = N_{\geq k}(G_n(1/n))$. We shall show that

(12) $$\mathbb{E}((N_{\geq k}/n)^2) \leq (\mathbb{E}(N_{\geq k}/n))^2 + o(1)$$

as $n \to \infty$. This will imply that $N_{\geq k}/n$ has variance $o(1)$, and hence that $N_k(G_n(1/n)) = N_{\geq k} - N_{\geq k+1}$ is concentrated about its mean.

Writing $c(v)$ for the number of vertices in the component of $G_n(1/n)$ containing a given vertex $v$, and letting $v$ and $w$ be independent random vertices of $G_n(1/n)$, (12) is equivalent to

(13) $$\mathbb{P}(c(v) \geq k, c(w) \geq k) \leq \mathbb{P}(c(v) \geq k)\mathbb{P}(c(w) \geq k) + o(1).$$

But this is more or less immediate from the fact that $\mathbb{P}(d(v,w) \leq 2k) = o(1)$, where $d(v,w)$ denotes graph distance in $G_n(1/n)$. Indeed, let us first fix $v$ and $w$. If $c(v) \geq k$, $c(w) \geq k$ and $d(v,w) > 2k$, then we can find disjoint sets of edges $E_v, E_w \subset G_n(1/n)$ such that the presence of all edges of $E_v$ in $G_n(1/n)$ is sufficient to guarantee that $c(v) \geq k$, and similarly with $v$ replaced by $w$. [In fact, if $d(v,w) > 2k$, then minimal witnesses for the events $c(v) \geq k$ and $c(w) \geq k$ must be disjoint.] In other words, writing, as usual, $A \square B$ for the event that two (increasing) events $A$ and $B$ have disjoint witnesses, if $d(v,w) > 2k$, then whenever the events $A = \{c(v) \geq k\}$ and $B = \{c(w) \geq k\}$ hold, so does $A \square B$. Hence

$$\mathbb{P}(A \cap B) \leq \mathbb{P}(d(v,w) \leq 2k) + \mathbb{P}(A \square B, d(v,w) > 2k)$$
$$\leq \mathbb{P}(d(v,w) \leq 2k) + \mathbb{P}(A \square B).$$

By the van den Berg–Kesten inequality [4] (for a more general inequality, see Reimer [20]), $\mathbb{P}(A \square B)$ is at most $\mathbb{P}(A)\mathbb{P}(B)$, so

$$\mathbb{P}(c(v) \geq k, c(w) \geq k) \leq \mathbb{P}(c(v) \geq k)\mathbb{P}(c(w) \geq k) + \mathbb{P}(d(v,w) \leq 2k).$$

Choosing $v, w$ uniformly at random and using the fact that $\mathbb{E}_{vw}\mathbb{P}(d(v,w) \leq 2k) = o(1)$, we obtain the bound (13) and hence (12), completing the proof of Lemma 4. □

Lemma 4 has an immediate corollary showing that the "right" number of vertices are in "small" components, as long as "small" is defined suitably.

COROLLARY 6. *Let $W : [0,1]^2 \to [0, \infty)$ be a graphon, and let $(G_n) = (\beta_{ij}(n))_{i,j \in [n]}$ be a sequence of weighted graphs with $G_n \to W$ and $\beta_{\max} = \sup_{i,j,n} \beta_{ij}(n)$ finite. Then*

$$\frac{1}{n} N_{>\omega}(G_n(1/n)) \xrightarrow{\mathrm{P}} \rho(W)$$

*as $n \to \infty$ whenever $\omega = \omega(n)$ tends to infinity sufficiently slowly.*

PROOF. By Lemma 4 we have $|N_k(G_n(1/n))/n - \mathbb{P}(|\mathfrak{X}_W| = k)| \xrightarrow{\mathrm{P}} 0$ for every fixed $k$. It follows that there is some function $\omega = \omega(n)$ with $\omega(n) \to \infty$ such that $\sum_{k=1}^{\omega} |N_k(G_n(1/n))/n - \mathbb{P}(|\mathfrak{X}_W| = k)| \xrightarrow{\mathrm{P}} 0$. Reducing $\omega$ if necessary, we may and shall assume that $\omega(n) = o(n)$. Let us say that a component is *small* if it has size at most $\omega(n)$ and *large* otherwise. Note that the number $N_{\leq \omega}(G_n(1/n))$ of vertices in small components satisfies

$$\frac{1}{n} N_{\leq \omega}(G_n(1/n)) = \sum_{k=1}^{\omega} \mathbb{P}(|\mathfrak{X}_W| = k) + o_{\mathrm{p}}(1) = \mathbb{P}(|\mathfrak{X}_W| \leq \omega(n)) + o_{\mathrm{p}}(1)$$
$$= \mathbb{P}(|\mathfrak{X}_W| < \infty) + o_{\mathrm{p}}(1).$$

(For the last step, recall that $\mathfrak{X}_W$ is defined without reference to $n$.) Hence, the number of vertices in large components is asymptotically $n\rho(W) = n\mathbb{P}(|\mathfrak{X}_W| = \infty)$, as claimed. □

In the case where $\|T_W\| \leq 1$, Theorem 1 follows from Corollary 6: we have $\rho(W) = 0$, so there are $o_{\mathrm{p}}(n)$ vertices in large components, and the largest component has size $o_{\mathrm{p}}(n)$. When $\|T_W\| > 1$, it remains to show that almost all vertices in large components are in a single giant component. For this we shall use a sprinkling argument.

**3. Joining up the large components.** In the light of Corollary 6, to prove the heart of Theorem 1, namely part (c), it remains to show that when $W$ is irreducible, almost all vertices in "large" components are in fact in a single giant component. Before doing this, we shall show that part (b) of Theorem 1, concerning the reducible case, follows from part (c). As the reducible case is rather uninteresting, we shall only outline the argument, omitting the straightforward details.

If $W$ is reducible, it is easy to check that $W$ may be decomposed into a finite or countably infinite sequence $W_1, W_2, \ldots$ of irreducible graphons in a natural sense. (For a formal statement and the simple proof, see Lemma 5.17 of [6].) Here each graphon $W_i$ is defined on $A_i \times A_i$, rather than on $[0,1]^2$, for some disjoint sets $A_i \subset [0,1]$. If $\|T_W\| > 1$, then there is an $i$ for which $\|T_{W_i}\| > 1$. In fact, since $W$ is bounded (by $\beta_{\max} < \infty$), only finitely many of the $W_j$ may have $\|T_{W_j}\| > 1$, so there is an $i$ with $\|T_{W_i}\| = \|T_W\|$. As noted in the Introduction, if $G_n \to W$ it is easy to check that we may find induced subgraphs $H_n$ of $G_n$ with $|H_n| = (\mu(A_i) + o(1))n$ and $H_n \to W_i$. Since $W_i$ is irreducible, assuming the irreducible case of Theorem 1, whp the graph $H_n(c/n)$, and hence $G_n(c/n)$, will contain a component with $\Theta(n)$ vertices.

To obtain the explicit constant in Theorem 1(b), it remains only to show that if $W$ (here $W_i$, rescaled to $[0,1]^2$) is an irreducible graphon with $\|T_W\| > 1$, then

$$\rho(W) \geq \frac{\|T_W\| - 1}{\|W\|_\infty}. \tag{14}$$

This crude bound is implicit in the results in [6]: indeed, $T_W$ has a positive eigenfunction $\psi$ (see [6], Lemma 5.15) with eigenvalue $\lambda = \|T_W\|$. Since

$$\|T_W\|\psi(x) = \lambda\psi(x) = (T_W\psi)(x) = \int_y W(x,y)\psi(y)\,dy \leq \|W\|_\infty \|\psi\|_1,$$

we have $\|T_W\|\|\psi\|_\infty \leq \|W\|_\infty \|\psi\|_1$. From [6], Remark 5.14, we have

$$\rho(W) \geq \frac{\|T_W\| - 1}{\|T_W\|} \frac{\|\psi\|_1}{\|\psi\|_\infty},$$

which then implies (14).

In the light of the comments above, from now on we assume that $W$ is irreducible. In joining up the large components to form a single giant component, we must somehow make use of this irreducibility. By an $(a,b)$-*cut* in $W$, we shall mean a partition $(A, A^c)$ of $[0,1]$ with $a \leq \mu(A) \leq 1 - a$ such that $\int_{A \times A^c} W \leq b$. We start with a simple lemma showing that irreducibility [no $(a,0)$-cut for any $0 < a \leq \frac{1}{2}$] implies an apparently stronger statement. Recall that our graphons are bounded by definition.

LEMMA 7. *Let $W$ be an irreducible graphon, and let $0 < a < \frac{1}{2}$ be given. There is some $b = b(W, a) > 0$ such that $W$ has no $(a, b)$-cut.*

PROOF. Define a measure $\nu$ on $X^2 = [0,1]^2$ by setting $\nu(U) = \int_U W(x, y)\,dx\,dy$ for each (Lebesgue-)measurable set $U \subset X^2$. Renormalizing, we may and shall assume that $W(x,y) \le 1$ for every $(x,y) \in X^2$, so that $\nu(B \times C) \le \mu(B)\mu(C)$. As $W$ is irreducible, we also have $w = \nu(X^2) > 0$.

Suppose that the assertion of the lemma is false. Then there is a sequence $(A_i, A_i^c)$ of pairs of complementary subsets of $X$ such that $a \le \mu(A_i) \le 1 - a$ and the $\nu$-measure of the cuts $C_i = (A_i \times A_i^c) \cup (A_i^c \times A_i)$ tends to 0. By selecting a subsequence, we may assume that $\nu(C_i) \le 2^{-i-1}w$ for every $i \ge 1$.

For $m \ge 1$, let $\mathcal{D}_m$ be the set of atoms of the partitions $\mathcal{P}_i = (A_i, A_i^c)$, $i = 1, \ldots, m$. Thus $D \in \mathcal{D}_m$ if and only if $D = \bigcap_{i=1}^m B_i$, with $B_i = A_i$ or $A_i^c$ for each $i$. Similarly, write $\mathcal{E}_n$ for the collection of atoms of the partitions $\mathcal{P}_n, \mathcal{P}_{n+1}, \ldots$. Since

$$X^2 = \bigcup_{i=1}^m C_i \cup \bigcup_{D \in \mathcal{D}_m} (D \times D)$$

we have

$$w = \nu(X^2) \le \sum_{i=1}^m \nu(C_i) + \sum_{D \in \mathcal{D}_m} \nu(D \times D) \le w/2 + \sum_{D \in \mathcal{D}_m} \mu(D)^2$$

$$\le w/2 + \max_{D \in \mathcal{D}_m} \mu(D) \sum_{D \in \mathcal{D}_m} \mu(D) = w/2 + \max_{D \in \mathcal{D}_m} \mu(D).$$

This shows that for each $m \ge 1$, we can find a $D_m \in \mathcal{D}_m$ with $\mu(D_m) \ge w/2$. Clearly, if $m < n$ and $D_n' \in \mathcal{D}_n$, then there is a (unique) $D_m' \in \mathcal{D}_m$ with $D_m' \supset D_n'$. Since each $\mathcal{D}_m$ is finite, by a standard compactness argument (repeated use of the pigeonhole principle) we may assume that $D_1 \supset D_2 \supset D_3 \supset \cdots$. Let $E_1 = \bigcap_{m=1}^\infty D_m$. Then $E_1 \in \mathcal{E}_1$ and $w/2 \le \mu(E_1) \le 1 - a$. For $n \ge 1$, let $E_n$ be the atom in $\mathcal{E}_n$ containing $E_1$; then $E_1 \subset E_2 \subset \cdots$ and $w/2 \le \mu(E_n) \le 1 - a$ for every $n$. Finally, set $E = \bigcup_{n=1}^\infty E_n$, so that $w/2 \le \mu(E) \le 1 - a$.

We claim that this set $E$ shows that $W$ is reducible. Indeed, for any $n$, $x \in E$ implies there is an $m \ge n$ with $x \in E_m$. Thus

$$E \times E^c \subset \bigcup_{m=n}^\infty (E_m \times E_m^c).$$

Since $E_m \times E_m^c \subset \bigcup_{i=m}^\infty C_i$, we have

$$\nu(E_m \times E_m^c) \le \sum_{i=m}^\infty \nu(C_i) \le \sum_{i=m}^\infty 2^{-i-1}w = 2^{-m}w,$$

so

$$\nu(E \times E^{\mathrm{c}}) \leq \inf_n \sum_{m=n}^{\infty} \nu(E_m \times E_m^{\mathrm{c}}) \leq \inf_n 2^{-n+1} w = 0,$$

contradicting our assumption that $W$ is irreducible. □

We are now ready to complete the proof of Theorem 1.

PROOF OF THEOREM 1. As before, we normalize so that $c = 1$. As noted at the start of the section, we may assume that $W$ is irreducible.

By Corollary 6, there is some $\omega = \omega(n)$ with $\omega(n) \to \infty$ and $\omega(n) = o(n)$ such that

$$\frac{1}{n} N_{>\omega}(G_n(1/n)) \xrightarrow{\mathrm{P}} \rho(W) = \mathbb{P}(|\mathfrak{X}_W| = \infty).$$

In particular, since the size $C_1$ of the largest component of $G_n(1/n)$ is at most the maximum of $\omega$ and $N_{>\omega}(G_n(1/n))$, for any $\varepsilon > 0$ we have $C_1 \leq (\rho(W) + \varepsilon)n$ whp. Similarly, the sum of the sizes of the two largest components is at most $(\rho(W) + 2\varepsilon)n$ whp. Since $\rho(W) = 0$ if $\|T_W\| \leq 1$, it remains only to show that, if $\|T_W\| > 1$ and $W$ is irreducible, then $C_1 \geq (\rho(W) - \varepsilon)n$ holds whp.

Theorem 6.4 of [6] which, like all results in Sections 5 and 6 of that paper, applies to all graphons (rather than the more restrictive kernels considered elsewhere in [6]), tells us that if $(W_k)$ is a sequence of graphons with $W_k$ tending up to $W$ pointwise, then $\rho(W_k) \to \rho(W)$. In particular, we have $\rho((1-\delta)W) \to \rho(W)$ as $\delta \to 0$. (If we are only interested in the existence of a giant component, rather than its size, then we may use instead Theorem 6.7 in [15].) It thus suffices to show that for any $\delta > 0$, whp we have

$$C_1/n \geq \rho((1-\delta)W) - 3\delta,$$

say. In doing so we may of course assume that $\delta$ is small enough that

(15) $$\rho((1-\delta)W) > 10\delta,$$

say. We shall also assume that $\delta \leq 1/100$.

Let $G' = G_n((1-\delta)/n)$ be the (unweighted) graph on $[n]$ in which the edges are present independently, and the edge $ij$ is present with probability $(1-\delta)\beta_{ij}/n$, where $\beta_{ij} = \beta_{ij}(n)$ is the weight of $ij$ in $G_n$. We may regard $G'$ as $G'_n(1/n)$, where $G'_n$ is obtained from $G_n$ by multiplying all edge weights by $1-\delta$. Since $G'_n \to (1-\delta)W$, by Corollary 6 there is some $\omega(n) \to \infty$ with $\omega(n) = o(n)$ such that

(16) $$N_{>\omega}(G')/n \geq \rho((1-\delta)W) - \delta$$

holds whp.

By an $(a,b)$-*cut* in an $n$-vertex weighted graph $G$ we shall mean a partition of the vertex set of $G$ into two sets $X$, $X^{\mathrm{c}}$ of at least $an$ vertices such that the total weight of edges from $X$ to $X^{\mathrm{c}}$ is at most $bn^2$. By Lemma 7, there is some $b > 0$ such that $W$ has no $(\delta, 2b)$-cut. We may and shall assume that $b < 1/10$, say. Since $G_n \to W$ in the cut metric [see (2)], it follows that if $n$ is large enough, which we shall assume from now on, then $G_n$ has no $(\delta, b)$-cut.

The graph $G_n(1/n)$ may be obtained from $G'$ by adding each nonedge $ij$ with a certain probability $s_{ij}$, independently of the other nonedges, where $s_{ij} \geq \delta \beta_{ij}/n$, and the inequality is strict unless $\beta_{ij} = 0$. In fact, we shall add each nonedge $ij$ with probability exactly $\delta \beta_{ij}/n$; the correction from this value to the true value of $s_{ij}$ only works in our favor.

Let us condition on $G'$, assuming as we may that (16) holds. Let $C_1, \ldots, C_r$ list the "large" components of $G'$, that is, the components with more than $\omega$ vertices. To complete the proof of Theorem 1, it suffices to show that, whp, all but at most $2\delta n$ vertices of $\bigcup C_i$ lie in a single component of $G_n(1/n)$. Let $B$ be the "bad" event that this does not happen, so we must show that $\mathbb{P}(B) = o(1)$. Since $G'$ is a subgraph of $G_n(1/n)$, whenever $B$ holds there is a partition $X \cup Y$ of $\{1, 2, \ldots, r\}$ such that $G_n(1/n)$ contains no path from $C_X = \bigcup_{i \in X} C_i$ to $C_Y = \bigcup_{i \in Y} C_i$, with $|C_X|, |C_Y| \geq 2\delta n$.

Given a weighted graph $G$, with edge weights $\beta_{vw}$, for $W \subset V(G)$ and $v \in V(G)$, we write

$$e(v, W) = e_G(v, W) = \sum_{w \in W} \beta_{vw}.$$

Similarly, for $V, W \subset V(G)$,

$$e(V, W) = e_G(V, W) = \sum_{v \in V} \sum_{w \in W} \beta_{vw}.$$

Unless stated otherwise, the quantities $e(v, W)$ and $e(V, W)$ will refer to the weighted graph $G = G_n$.

Fix $G'$ (and hence $C_1, \ldots, C_r$) and a partition $X, Y$ of $\{1, 2, \ldots, r\}$ for which $|C_X|, |C_Y| \geq 2\delta n$. We shall inductively define an increasing sequence $S_0, S_1, \ldots, S_\ell$ of sets of vertices of $G_n$, in a way that depends on $C_X$ and on $G_n$, but not on the "sprinkled" edges of $G_n(1/n) \setminus G'$. We start with $S_0 = C_X$, noting that $|S_0| \geq 2\delta n$. We shall stop the sequence when $|S_t|$ first exceeds $(1-\delta)n$. Thus, in defining $S_{t+1}$ from $S_t$, we may assume that $\delta n \leq |S_t| \leq (1-\delta)n$. Since $G_n$ has no $(\delta, b)$-cut, we have

$$\sum_{v \notin S_t} e(v, S_t) = e(S_t^{\mathrm{c}}, S_t) \geq bn^2.$$

Let

$$T_{t+1} = \{v \notin S_t : e(v, S_t) \geq bn/2\}.$$

PERCOLATION ON DENSE GRAPH SEQUENCES 21

As all edge weights in $G_n$ are bounded by $\beta_{\max}$, we have $e(v, S_t) \leq \beta_{\max}|S_t| \leq \beta_{\max} n$ for any $v$, so

$$bn^2 \leq e(S_t^c, S_t) \leq \frac{bn}{2}|V(G_n) \setminus (S_t \cup T_{t+1})| + \beta_{\max} n|T_{t+1}| \leq \frac{bn^2}{2} + \beta_{\max} n|T_{t+1}|.$$

Hence, $|T_{t+1}| \geq \frac{bn}{2\beta_{\max}}$. Set $S_{t+1} = S_t \cup T_{t+1}$, and continue the construction until we reach an $S_\ell$ with $|S_\ell| \geq (1-\delta)n$. Note that $\ell \leq \frac{2\beta_{\max}}{b} = O(1)$.

We shall now uncover the "sprinkled" edges between $T_t$ and $S_{t-1}$, working backwards from $T_\ell$. It will be convenient to set $T_0 = S_0$, so $S_t = \bigcup_{t'=0}^t T_{t'}$. Since $|S_\ell| \geq (1-\delta)n$, while $|C_Y| \geq 2\delta n$, the set $S_\ell$ contains at least $\delta n$ vertices from $C_Y$. Since $S_0 = T_0 = C_X$ is disjoint from $C_Y$, it follows that there is some $t_0$, $1 \leq t_0 \leq \ell$, for which $T_{t_0}$ contains a set $Y_0$ of at least $\delta n/\ell$ vertices of $C_Y$. Passing to a subset, we may assume that

$$|Y_0| = \min\left\{\frac{\delta n}{\ell}, \frac{bn}{10\beta_{\max}}\right\} + O(1),$$

so $|Y_0| = \Theta(n)$ but $|Y_0| \leq bn/(10\beta_{\max}) \leq |T_t|/5$ for $1 \leq t \leq \ell$.

Next, we construct a set $X_0 \subset S_{t_0-1}$ with $|X_0| \geq \delta b|Y_0|/5$ such that every $x \in X_0$ is joined to some $y \in Y_0$ by an edge of $G_n(1/n) \setminus G'$. We start with the observation that for every vertex $y \in Y_0$ we have $e(y, S_{t_0-1}) \geq bn/2$, since $y \in T_{t_0}$. Hence the expected number of edges of $G_n(1/n) \setminus G'$ from $y$ to $S_{t_0-1}$ is at least $\delta b/2$, and the probability that at least one such edge is present is at least $1 - \exp(-\delta b/2) \geq \delta b/4$. Furthermore, this conclusion remains true (with $\delta b/2$ replaced by $\delta b/3$) even if we exclude a subset of $S_{t_0-1}$ of size at most $|Y_0|$, corresponding to at most one neighbor $x' \in S_{t_0-1}$ of each vertex $y' \in Y_0$ previously considered. [To see this, we note that for every $\tilde{X}_0 \subset S_{t_0-1}$ with $|\tilde{X}_0| \leq |Y_0|$, we have $e(y, S_{t_0-1} \setminus \tilde{X}_0) \geq bn/2 - \beta_{\max}|Y_0| \geq bn/2 - bn/10 \geq bn/3$.] Using independence of edges from different vertices $y$, and the concentration of the binomial distribution, it follows that with probability at least $1 - \exp(-\Theta(n))$, we find a set $X_0$ of at least $\delta b|Y_0|/5$ vertices of $S_{t_0-1}$ such that every $x \in X_0$ is joined to some $y \in Y_0$ by an edge of $G_n(1/n) \setminus G'$.

Since $|X_0| \geq \delta b|Y_0|/5$, there is some $t_1 < t_0$ such that $Y_1 = X_0 \cap T_{t_1}$ contains at least $\delta b|Y_0|/(5\ell)$ vertices. If $t_1 \geq 1$ then, arguing as above, with probability $1 - \exp(-\Theta(n))$ we find a $t_2$ and a set $Y_2$ of at least $\delta^2 b^2|Y_0|/(5\ell)^2$ vertices of $T_{t_2}$ joined in $G_n(1/n)$ to $Y_1$, and so on. As the sequence $t_0, t_1, \ldots$ is decreasing, this process terminates after $s \leq \ell$ steps with $t_s = 0$. Hence, with probability $1 - \exp(-\Theta(n))$ we find a set $Y_s$ of at least $(\delta b/(5\ell))^\ell |Y_0| = \Theta(n) > 1$ vertices of $T_0 = S_0 = C_X$ joined in $G_n(1/n)$ by paths to vertices in $C_Y$. In particular, the probability that there is no path in $G_n(1/n)$ from $C_X$ to $C_Y$ is exponentially small.

Recalling that $r \leq n/\omega = o(n)$, the number of possible partitions $X, Y$ of the components $C_1, \ldots, C_r$ is at most $2^r = \exp(o(n))$, so the probability of the bad event $B$ is $o(1)$, as required. $\square$

**4. Stronger bounds on the small components.** In this section we prove Theorem 3, considering the subcritical and supercritical cases separately.

4.1. *The subcritical case.* We start by proving the first statement of Theorem 3, restated below for ease of reference.

THEOREM 8. *Let $(G_n) = (\beta_{ij}(n))_{i,j \in [n]}$ be a sequence of weighted graphs with $|G_n| = n$ and $\beta_{\max} = \sup_{i,j,n} \beta_{ij}(n) < \infty$ converging to a graphon $W$. If $\|T_W\| < 1$, then there is a constant $A$ such that $C_1(G_n(1/n)) \leq A \log n$ holds whp.*

PROOF. Let $\delta$ be a positive constant chosen so that $(1+\delta)\|T_W\| < 1$. Let $W_n = W_{G_n}$ be the piecewise constant graphon naturally associated to $G_n$, and let $W'_n = (1+\delta)G_n$. We claim that if $n$ is large enough, then the neighborhood exploration processes associated to a random vertex of $G_n$ may be coupled with $\mathfrak{X}_{W'_n}$, viewed as an $n$-type branching process, so that the latter dominates.

In the exploration process, we start with a random vertex of $G_n$. Having reached a vertex $i$, we check for possible "new" neighbors of $i$ not yet reached from other vertices. The chance that $j$ is such a new neighbor is either $\beta_{ij}(n)/n$ or 0, depending on whether $j$ has been previously reached or not. In particular, this process is dominated by (may be regarded as a subset of) a binomial $n$-type process in which we start with a particle of a random type, and each particle of type $i$ has a Bernoulli $B(\beta_{ij}(n)/n)$ number of children of type $j$, independently of everything else. The process $\mathfrak{X}_{W'_n}$ may be described in exactly the same terms except that the number of children of type $j$ has a $\mathrm{Po}((1+\delta)\beta_{ij}(n)/n)$ distribution. As $\beta_{ij}(n)$ is uniformly bounded and $\delta$ is fixed, this Poisson distribution dominates the corresponding Bernoulli distribution for all large enough $n$.

Although the branching processes $\mathfrak{X}_{W'_n}$ have different kernels, these kernels are uniformly bounded. Furthermore, since $W$ and the $W_n$ are uniformly bounded and $W_n \to W$ in the cut norm, it is easy to check (for example, by considering step function approximations to eigenfunctions of the compact operators $T_{W_n}$ and $T_W$) that $\|T_{W_n}\| \to \|T_W\|$. (In fact, as mentioned in the Introduction, much more is true—in [15] it is shown that the normalized spectra converge.) Hence, for all sufficiently large $n$, the norms $\|T_{W'_n}\|$ are bounded by some constant strictly less than 1. It is a standard result that the branching processes $\mathfrak{X}_{W'_n}$, associated to uniformly bounded, uniformly subcritical kernels $W'_n$ exhibit uniformly exponential decay; in other words,

$$\mathbb{P}(|\mathfrak{X}_{W'_n}| \geq k) \leq \exp(-ak) \tag{17}$$

for all sufficiently large $n$ and all $k \geq 1$, where $a > 0$ is constant. [This can be shown by considering $\mathbb{E}(e^{t|\mathfrak{X}_{W'_n}|})$.]

Finally, from the coupling above, the expected number of vertices of $G_n(1/n)$ in components of size at least $A \log n$ is at most $n \mathbb{P}(|\mathfrak{X}_{W'_n}| \geq A \log n)$; for $A = 2/a$, say, (17) tells us that this is at most $n \exp(-a\frac{2}{a} \log n) = 1/n = o(1)$, so whp there are no such vertices. □

4.2. *The supercritical case.* In the supercritical case we shall show that there is a constant $A = A(W)$ such that whp the second largest component of $G_n(1/n)$ contains at most $A \log n$ vertices. More precisely, we shall prove the second part of Theorem 3, that is, the following result.

THEOREM 9. *Let $(G_n) = (\beta_{ij}(n))_{i,j \in [n]}$ be a sequence of weighted graphs with $|G_n| = n$ and $\beta_{\max} = \sup_{i,j,n} \beta_{ij}(n) < \infty$ converging to a graphon $W$. Suppose that $\|T_W\| > 1$ and that $W$ is irreducible. Then there is a constant $A$ such that $C_2(G_n(1/n)) \leq A \log n$ holds whp, where $C_2(G)$ denotes the number of vertices in the second largest component of a graph $G$.*

The basic strategy of the proof is to use an idea from [6], although there will be considerable difficulties in adapting it to the present context. Let us first give a rough description of this idea. Note that we have already shown in Theorem 1 that $G_n(1/n)$ has whp a unique component of order $\Theta(n)$, the "giant" component. All other components are "small," that is, of order $o_p(n)$.

Suppose that we have a "supercritical" random graph $H$ on $n$ vertices [here $H = G_n(1/n)$], and let $k$ be a large constant to be chosen later. Let us select $n/k$ of the vertices of $H$ at random to be "left" vertices, the remaining vertices being "right" vertices; we do this *before* deciding which edges are present in the random graph $H$. If $H$ has probability $\varepsilon$ of containing a small component with at least $Ak \log n$ vertices, then [considering a random partition of $V(H)$ into $k$ parts] with probability at least $\varepsilon/k$ the graph $H$ contains a small component that in turn contains at least $A \log n$ left vertices. Thus, it suffices to show that whp any component of $H$ containing at least $A \log n$ left vertices is the unique component of $H$ with size $\Theta(n)$.

If $k$ is chosen large enough, then the subgraph induced by the right vertices already contains a component of size $\Theta(n)$. Uncovering the subgraphs $H_R$ and $H_L$ of $H$ induced by the right and left vertices, respectively, and all edges between the small components of $H_R$ and the left vertices, we have *already revealed* the small components of $H$: writing $H'$ for the subgraph of $H$ formed by the edges revealed so far, any small component of $H$ is a small component of $H'$. All that remains is to uncover the edges of $H$ between the unique giant component of $H_R$ and the left vertices; adding these edges to $H'$ will cause certain components of $H'$ to merge into the giant component but have no other effect. If the edge probabilities in $H$ are bounded below by $c/n$, $c > 0$, and $A$ is chosen to be large enough, then it is very unlikely that any component of $H'$ with $A \log n$ left vertices fails to merge into the giant component, so it is very unlikely that $H$ has a small component with at least $A \log n$ left vertices.

This argument can be applied as it is to $H = G_n(1/n)$ if the edge weights $\beta_{ij}(n)$ in $G_n$ are bounded away from zero. However, this is typically not the case. Indeed, the main interest is when $G_n$ is a graph rather than a weighted graph, so many edge weights will be zero. To overcome this difficulty we shall use regularity instead of a lower bound on individual edge probabilities.

Our notation for regularity is standard. Thus, for disjoint sets $A$, $B$ of vertices of a weighted graph $G$, we write $e(A,B) = e_G(A,B)$ for the total edge weight from $A$ to $B$ in $G$. Also, $d(A,B) = e(A,B)/(|A||B|)$ is the *density* of the pair $(A,B)$. An *$\varepsilon$-regular pair* is a pair $(A,B)$ of sets of vertices of $G$ such that

$$|d(A',B') - d(A,B)| \leq \varepsilon$$

whenever $A' \subset A$ and $B' \subset B$ with $|A'| \geq \varepsilon|A|$ and $|B'| \geq \varepsilon|B|$. A partition $\mathcal{P}$ of $V(G)$ into sets $P_1, \ldots, P_M$ is *$\varepsilon$-regular* if $|P_i| = \lfloor n/M \rfloor$ or $\lceil n/M \rceil$ for every $i$, $|P_i| \leq \varepsilon n$, and all but at most $\varepsilon M^2$ of the pairs $(P_i, P_j)$, $i \neq j$, are $\varepsilon$-regular. Szemerédi's Lemma [21] tells us that, for any $\varepsilon > 0$, there exist $M = M(\varepsilon)$ and $n_0 = n_0(\varepsilon)$ such that any graph $G$ with $n \geq n_0$ vertices has an $\varepsilon$-regular partition $\mathcal{P}$ into exactly $M$ classes. (Most formulations give at most $M$ classes; applying this weaker form with a slightly smaller $\varepsilon$, and then randomly partitioning classes into smaller classes of the desired size, it is easy to deduce the stated form.) Note that while Szemerédi's original lemma applies to unweighted graphs, the extension to graphs with bounded edge weights is immediate (either by adapting the proof, or by first replacing all edge weights with multiples of $\varepsilon/100$, say).

Given a Szemerédi partition as above, we write $G/\mathcal{P}$ for the weighted graph with vertex set $P_1, \ldots, P_M$ in which the weight of the edge $P_i P_j$ is $d(P_i, P_j)$ if $(P_i, P_j)$ is $\varepsilon$-regular, and zero otherwise. If the edge weights in the graph $G$ are bounded by $\beta_{\max} < \infty$, then it is easy to check that the cut metric distance $\delta_\square(G, G/\mathcal{P})$ is at most $4\varepsilon(\beta_{\max} + 1)$, say.

For the rest of the paper we fix $G_n$ with $G_n \to W$, $\beta_{\max} < \infty$, $W$ irreducible and $\|T_W\| > 1$. Note that $\|T_W\| \leq \beta_{\max}$, so $\beta_{\max} > 1$. We also fix a constant $k \geq 2$ with $(1 - 1/k)^2 \|T_W\| > 1$. Let $L$ be a set of $n/k$ vertices of $G = G_n$ chosen uniformly from among all such sets, and set $R = V(G) \setminus L$: these are the sets of *left* and *right* vertices. Here and throughout we ignore rounding to integers, which clutters the notation without affecting the proofs.

In what follows we shall consider several different graphs defined in terms of $G$, $L$ and $R$, some random and some not. We write $G_L$ and $G_R$ for the subgraphs of $G = G_n$ induced by $L$ and $R$, respectively. We use a superscript minus to denote graphs in which the edge weights have been multiplied by $(1 - 1/k)$ (i.e., reduced slightly), and we replace $G$ by $H$ to denote the random "subgraph" of a weighted graph ($G$, $G_R^-$, etc.) obtained by selecting each edge $ij$ with probability given by its weight divided by $n$. In particular, we shall consider the following graphs at various stages:

$$G = G_n,$$
$$H = G_n(1/n),$$
$$G_R = G_n[R],$$
$$H_R = H[R] = G_n(1/n)[R],$$
$$H_L = H[L] = G_n(1/n)[L],$$
$$G_R^- = (1 - 1/k)G_n[R],$$
$$H_R^- = (1 - 1/k)G_n[R](1/n) = G_n[R]((1 - 1/k)/n)$$

as well as

$$C_1 = \text{largest component of } H_R$$

and

$$C_1^- = \text{largest component of } H_R^-.$$

We shall also consider the graph $H'$ defined as in the outline proof above: $H'$ will be the subgraph of $H$ consisting of $H_R \cup H_L$ together with all edges of $H$ joining $L$ to $R \setminus C_1$.

LEMMA 10. *Let $G_n$, $W$, $k \geq 2$, $L$ and $R$ be as above. Then, with probability 1, $G_R \to W$. Furthermore, for any $\delta > 0$ the largest component $C_1^-$ of $H_R^-$ satisfies*

$$|C_1^-|/n \geq (1 - 1/k)\rho((1 - 1/k)^2 W) - \delta$$

*whp.*

PROOF. The first statement is a more or less immediate consequence of Szemerédi's regularity lemma; indeed, using this lemma and a simple sampling argument, one easily shows that $\delta_\square(G, G_R) \to 0$ in probability (see Theorem 2.9 of [14] for an explicit bound on the error term). Although convergence in probability is all we need here, the error probability can be made small enough to ensure convergence with probability 1. Since $G \to W$, we then have $G_R \to W$.

The second part follows by Theorem 1, noting that $G_R^-$ has $n' = (1-1/k)n$ vertices, so $H_R^- = G_R^-(1/n) = G_R^-((1-1/k)/n')$, and that $G_R \to W$ trivially implies $G_R^- \to (1-1/k)W$. □

We next note that irreducibility of $W$ implies that $G$ cannot have too many "low-degree" vertices. The constants in this lemma are not written in the simplest form, but rather in the form that we shall use later.

LEMMA 11. *Let $G_n \to W$ with $W$ irreducible and $\beta_{\max} = \sup_{i,j,n} \beta_{ij}(n) < \infty$, and let $k \geq 2$ be constant. There is a constant $\sigma > 0$ such that, for all large enough $n$, at most $n/(10\beta_{\max})$ vertices of $G_n$ have weighted degree less than $50k\beta_{\max}\sigma n$.*

PROOF. Set $a = 1/(10\beta_{\max})$. By Lemma 7 there is a $\sigma$ such that $W$ has no $(a, 6k\sigma)$-cut. Since $G_n \to W$, it follows that, for large enough $n$, the graph $G_n$ has no $(a, 5k\sigma)$-cut. But then the conclusion of the lemma follows: otherwise, let $A$ be a set of exactly $n/(10\beta_{\max})$ vertices of $G_n$ with $d_v \leq 50k\beta_{\max}\sigma n$ for every $v \in A$. Then $e(A, A^c) \leq 50k\beta_{\max}\sigma n |A| = 5k\sigma n^2$, so $(A, A^c)$ is an $(a, 5k\sigma)$-cut in $G_n$. □

From now on, we fix a constant $\sigma > 0$ for which Lemma 11 holds. We assume, as we may, that

$$\sigma < 10^{-10}/\beta_{\max}^3 \tag{18}$$

and that

$$\sigma < (1 - 1/k)\rho((1-1/k)^2 W)/3. \tag{19}$$

Our next lemma shows that when we split $G = G_n$ into left and right vertices, most vertices on the left will send a reasonable weight of edges to $C_1$, the giant component of the random subgraph on the right, that is, the largest component of $H_R = G_R(1/n) = G_n(1/n)[R]$. By Lemma 10 and our assumption (19), we have

$$|C_1| \geq |C_1^-| \geq 2\sigma n$$

whp. Also, by Theorem 1, $C_1$ is whp the unique component of $H_R$ with $\Theta(n)$ vertices.

LEMMA 12. *Under the assumptions above there is a constant $\gamma > 0$ such that whp the number of vertices $v \in L$ with $e_G(v, C_1) \leq \gamma n$ is at most $3\sigma n$.*

PROOF. The proof of this technical lemma is somewhat involved. The basic idea is to take a Szemerédi partition and use regularity. The problem is that a priori we have no control over the distribution of $C_1$ with respect to this partition: it could sit almost entirely within a few parts, and there might be many other parts with low density to these parts. The solution is to start from $C_1^-$ and use irreducibility to show that the extra sprinkled edges expand $C_1^-$ to a set $(C_1)$ that has a decent number of vertices in almost all parts of the partition.

By Lemma 7, there is a constant $b > 0$ such that $W$ has no $(\sigma, 3b)$-cut. We may and shall assume that $b < 1/1000$, say. Recalling that $\beta_{\max} > 1$, set

$$\varepsilon = \frac{\sigma}{2}\left(\frac{b^2}{33k\beta_{\max}}\right)^{\beta_{\max}/b}$$

and apply Szemerédi's lemma to find, for all sufficiently large $n$, an $\varepsilon$-regular partition $\mathcal{P} = \mathcal{P}_n$ of the weighted graph $G = G_n$ into $M$ classes $P_1, \ldots, P_M$. As the partition $\mathcal{P}$ depends on $G$ only, not on the random choice of $L$, whp each class $P_i$ satisfies $|P_i \cap L| = n/(kM) + o(n)$. Also, whp every vertex $v$ satisfies $e_G(v, L) = d_v/k + o(n)$. From now on we condition on $L$, assuming these two properties.

Let $G/\mathcal{P}$ be defined as above, noting that $\delta_\square(G/\mathcal{P}, W) \leq C\varepsilon$ where $C = 4(\beta_{\max} + 1)$ is constant. As $b < 1/1000$, we have $C\varepsilon < b$. Hence, if $n$ is large enough, which we assume from now on, the graph $G/\mathcal{P}$ has no $(\sigma, 2b)$-cut.

Recall that $H_R^-$ is the graph on $R$ in which each edge is present independently, and the probability that $ij$ is an edge is $(1 - 1/k)\beta_{ij}(n)/n$. Note that the subgraph $H_R = G(1/n)[R]$ of $H = G_n(1/n)$ induced by the vertices in $R$ may be obtained from $H_R^-$ by "sprinkling," that is, by adding each nonedge $ij$ of $H_R^-$ with probability (a little larger than) $\beta_{ij}/(kn)$. Let $C_1^-$ be the largest component of $H_R^-$. As noted above, by Lemma 10 and our choice of $\sigma$, we have $|C_1^-| \geq 2\sigma n$ whp. From now on we condition on $H_R^-$, assuming that this holds.

Let $S_0'$ be the set of classes $P_i$ with $|C_1^- \cap P_i| \geq \sigma n/M$. As

$$2\sigma n \leq |C_1^-| \leq (n/M)|S_0'| + (\sigma n/M)(M - |S_0'|) \leq (n/M)|S_0'| + \sigma n,$$

we have $|S_0'| \geq \sigma M$. Let $S_0$ be an arbitrary subset of $S_0'$ of size exactly $\sigma M$ (ignoring rounding, as usual). We shall inductively define an increasing sequence $S_0 \subset S_1 \subset \cdots \subset S_\ell$ of sets of classes of $\mathcal{P}$, stopping the first time we reach an $S_t$ with $|S_t| \geq (1 - \sigma)M$. In doing so, we shall write $T_i$ for $S_i \setminus S_{i-1}$.

Having defined $S_i$, with $\sigma M \leq |S_i| \leq (1 - \sigma)M$, let $T_{i+1}$ be the set of classes $P_j$ with $e(P_j, S_i) \geq bM$, where $e(\cdot, \cdot)$ counts the weight of edges *in the graph* $G/\mathcal{P}$. As $G/\mathcal{P}$ has no $(\sigma, 2b)$-cut, we have $e(S_i^c, S_i) \geq 2bM^2$. Since $e(P_j, S_i) \leq \beta_{\max}|S_i| \leq \beta_{\max}M$ for every $j$, it follows that $|T_{i+1}| \geq bM/\beta_{\max}$. Set $S_{i+1} = S_i \cup T_{i+1}$ and continue until we reach $S_\ell$ with $|S_\ell| \geq (1 - \sigma)M$. Since there are only $M$ classes in total, and $|T_i| \geq bM/\beta_{\max}$ for every $i$, the process just defined stops after $\ell \leq \beta_{\max}/b$ steps.

CLAIM 13. *Whp every class $P_j$ in $S_\ell$ contains at least $\varepsilon n/M$ vertices of $C_1$, the giant component of $H_R$.*

To show this, we shall use induction to prove the stronger statement that, whp, every class $P_j \in T_i$ contains at least $c_i n/M$ vertices of a certain set $C'_i \subset C_1$ that we shall define in a moment, where

$$c_i = \sigma\left(\frac{b^2}{33k\beta_{\max}}\right)^i$$

and it is convenient to set $T_0 = S_0$.

Set $C'_0 = C_1^-$, so the base case $i = 0$ holds by the definition of $S_0 = T_0$. In proving the induction step, we shall use the "sprinkled" edges between $T_{i+1} \cap R$ and $S_i \cap R$; we define $C'_{i+1}$ to be the set obtained from $C'_i$ by adding all vertices of $T_{i+1}$ joined directly to $C'_i$ by such sprinkled edges.

For the induction step, let $P \in T_{i+1}$ be a class of the partition $\mathcal{P}$. By definition of $T_{i+1}$, we have $e(P, S_i) \geq bM$, where $e(\cdot,\cdot)$ counts the weight of edges in $G/\mathcal{P}$. Since $e(P_j, P_k) \leq \beta_{\max}$ for every $j$, $k$, it follows that there is a set $Q$ of least $bM/(2\beta_{\max})$ classes $P_j \in S_i$ with $e(P, P_j) \geq b/2$ for every $P_j \in Q$. By definition of $G/\mathcal{P}$, each pair $(P, P_j)$, $P_j \in Q$, is $\varepsilon$-regular (in $G_n$) with density at least $b/2$. Set $P' = \bigcup_{P_j \in Q} P_j$. From standard properties of $\varepsilon$-regularity, it follows that the pair $(P, P')$ is $(2\varepsilon)$-lower regular with density at least $b/4$. By the induction hypothesis, each $P_j \in S_i$ contains at least $c_i n/M = c_i |P_j|$ vertices of $C'_i$. Hence $|P' \cap C'_i| \geq c_i |P'| \geq 2\varepsilon|P'|$. By regularity, it follows that

$$(20) \qquad |\{v \in P : e_G(v, C'_i \cap P') \leq b|C'_i \cap P'|/5\}| \leq 2\varepsilon|P| = 2\varepsilon n/M.$$

In particular, of the $(1 + o(1))(1 - 1/k)n/M$ vertices of $P \cap R$, at least $N = n/(3M)$, say, have

$$(21) \quad \begin{aligned} e_G(v, C'_i) &\geq e_G(v, C'_i \cap P') \geq b|C'_i \cap P'|/5 \geq c_i b|P'|/5 \\ &\geq c_i b^2 n/(10\beta_{\max}), \end{aligned}$$

where for the last inequality we used $|Q| \geq bM/(2\beta_{\max})$.

A standard calculation using concentration of the Binomial distribution implies that, whp, at least $Nc_i b^2/(11k\beta_{\max}) = c_{i+1}|P|$ of these vertices $v$ are joined to $C'_i$ by a sprinkled edge, so $|C'_{i+1} \cap P| \geq c_{i+1}|P|$, completing the induction argument and hence (as $C'_i \subset C_1$ for all $i$) proving the claim.

The proof of Lemma 12 is also essentially complete. Set $\gamma = c_\ell b^2/(10\beta_{\max})$. Since $|S_0|$, $|S^c_\ell| \leq \sigma M$, it suffices to prove that, whp, for each $P \in T_j$, $1 \leq j \leq \ell$, there are at most $\sigma|P|$ vertices $v$ of $P \cap L$ with $e_G(v, C_1) \leq \gamma n$. Since $C'_i \subset C_1$ for every $i$, while $c_i \geq c_\ell$ and $2\varepsilon < \sigma$, this is immediate from (20) and (21). □

In proving Theorem 9, we shall first uncover all edges of $H = G_n(1/n)$ between vertices in $R$. In addition to revealing the giant component $C_1$ of $H_R = G_n(1/n)[R]$, this also reveals the small components of $H_R$. It will turn out that certain small components cause difficulty in the proof. Let us say that a vertex $v$ of $G = G_n$ has *low degree* if $e(v, L) \leq 49\beta_{\max}\sigma n$, where $\sigma$ is the constant in Lemma 11. We write $R^- \subset R$ for the set of low-degree vertices in $R$. Note that from Lemma 11 and the randomness of our partition, we have

$$|R^-| \leq n/(10\beta_{\max})$$

whp. Let us say that a component $C$ of $H_R$ is *annoying* if $|C \cap R^-| > |C|/2$.

LEMMA 14. *Let $A_s$ be the number of vertices of $H_R$ in annoying $s$-vertex components. Then whp we have $A_s \leq n \exp(-s/200)$ for every $s \geq 1$.*

PROOF. Throughout we confine our attention to the subgraph of $G$ induced by the vertices in $R$. We shall condition on $R$, assuming, as we may, that $|R^-| \leq n/(10\beta_{\max})$. For any $s$ we have $A_s \leq 2|R^-| \leq n/(5\beta_{\max}) \leq n/5$, so we may assume that $s \geq 200$, say.

Let $v$ be a random vertex of $R$, and let $C_v$ denote the component of $H_R$ containing $v$. We shall first show that, for some constant $a > 0$, we have

$$p_s = \mathbb{P}(|C_v| = s, |C_v \cap R^-| \geq |C_v|/2) \leq \exp(-as).$$

Note that

$$(22) \qquad p_s \leq \mathbb{P}(|C_v| \geq s, |C_v \cap R^+| \leq s/2),$$

where $R^+ = R \setminus R^-$. Let us explore the component $C_v$ in the usual way, writing $v_1, v_2, \ldots, v_t$ for the vertices of $C_v$ in the order we reach them, so $v_1 = v$. For technical reasons, we continue the sequence $(v_i)$ by starting a new exploration at a new random vertex whenever we exhaust the component currently being explored.

If $|C_v| \geq s$ and $|C_v \cap R^+| \leq s/2$, then at least $s/2 - 1$ of the vertices $v_2, \ldots, v_s$ are in $R^-$. In particular, at least $s/2 - 1$ of the children of $v_1, \ldots, v_s$ are in $R^-$. Let $\ell_t$ denote the number of children of $v_t$ that are in $R^-$. Then

$$p_s \leq \mathbb{P}\left(\sum_{t=1}^s \ell_t \geq s/2 - 1\right).$$

Now, as $|R^-| \leq n/(10\beta_{\max})$, for any vertex $w$ we have $e_G(w, R^-) \leq n/10$. In particular, when we test edges from a vertex $v_t$ to vertices not yet reached in the exploration, the chance of finding more than $k$ edges to $R^-$ in the random subgraph $H_R = G_R(1/n)$ of $G_R$ is at most $10^{-k}$. Hence,

$$p_s \le \sum_{r \ge s/2 - 1} \binom{s+r-1}{s-1} 10^{-r},$$

with the first factor coming from the number of sequences $k_1, \ldots, k_s$ with $k_i \ge 0$ and $\sum k_i = r$. A simple calculation shows that $p_s \le \exp(-s/100)$ for all $s \ge 200$, say.

Let $a = 1/100$, and fix an $s \ge 200$. Then $\mathbb{E}(A_s) = |R| p_s \le n \exp(-as)$, and it is easy to show that $\mathrm{Var}(A_s/n) = o(1)$. [This follows easily from the observation that the probability that two fixed vertices $v$ and $w$ are in the same component of size $s$ is $o(1)$, while for disjoint sets $X$ and $Y$ of $s$ vertices, the events that $X$ and $Y$ are the vertex sets of the components containing $v$ and $w$ are almost independent.] Setting $c = a/2 = 1/200$, it follows that $A_s \le n \exp(-cs)$ holds whp for each fixed $s$. Hence, there is some $s_0(n)$ tending to infinity such that whp the bound $A_s \le n \exp(-cs)$ holds simultaneously for all $s \le s_0$. For $s \ge s_0$ we simply use Markov's inequality to note that

$$\mathbb{P}(\exists s \ge s_0 : A_s \ge n \exp(-cs)) \le \sum_{s \ge s_0} \mathbb{P}(A_s \ge n \exp(-cs))$$
$$\le \sum_{s \ge s_0} \frac{n \exp(-as)}{n \exp(-cs)}$$
$$= \sum_{s \ge s_0} \exp(-cs) \to 0. \qquad \square$$

We are finally ready to prove Theorem 9.

PROOF OF THEOREM 9. Let $G_n \to W$ be a sequence of weighted graphs satisfying the assumptions of the theorem, and define $k$, $L$, $R$ and $\sigma$ as above. As before, we write $H = G_n(1/n)$ for the random graph whose component distribution we are studying.

As noted earlier, if $H$ has probability $\varepsilon$ of containing a small component with at least $Ak \log n$ vertices, then with probability at least $\varepsilon/k$, $H$ has such a component containing at least $A \log n$ vertices from $L$. By Theorem 1, $H$ has whp a unique component with $\Theta(n)$ vertices, so it suffices to prove that, for some constants $A$ and $\delta > 0$, whp every component of $H$ containing at least $A \log n$ vertices of $L$ has size at least $\delta n$; we shall prove this with $\delta = 2\sigma$ and $A$ a (large) constant to be chosen later.

As above, define $C_1$ to be the largest component of $H_R = H[R]$. Let $H'$ be the subgraph of $H$ consisting of all edges within $R$, all edges within $L$ and all edges between $L$ and vertices of $R \setminus C_1$. As $H$ is formed from $H'$ by adding some edges between $C_1$ and other components, whp every small component of $H$ is a component of $H'$. In particular, it suffices to prove that whp every component $C$ of $H'$ with at least $A \log n$ vertices in $L$ is joined to $C_1$ in $H$.

Let $\gamma$ be as in Lemma 12 above. For the rest of the proof we condition on $H_R$, assuming as we may that the property described in Lemma 14 holds. Call a vertex $v \in L$ *bad* if $e_G(v, C_1) \leq \gamma n$, and let $B \subset L$ be the set of bad vertices. By Lemma 12, we may assume that $|B| \leq 3\sigma n$.

It suffices to prove the following claim, which should be taken to hold conditional on $H_R$, under the assumptions above.

CLAIM 15. *Let $v$ be a random vertex of $L$, and let $C_v$ be the component of $H'$ containing $v$. There is a constant $c > 0$ such that the (conditional) probability that $|C_v \cap L| = s$ and $|C_v \cap B| \geq |C_v \cap L|/2$ is bounded above by $\exp(-cs)$ whenever $0 \leq s \leq \sigma n$.*

Indeed, to deduce Theorem 9 from Claim 15, note that if $|C_v \cap L| = s \geq A \log n$ and $|C_v \cap B| \leq |C_v \cap L|/2$, then when we uncover the (so far untested) edges between $L$ and $C_1$, the probability that there is no edge from $C_v$ to $C_1$ in $H$ is bounded by

$$\exp(-\gamma |C_v \cap (L \setminus B)|) \leq \exp(-s\gamma/2) \leq \exp(-A\gamma \log n/2),$$

which we can make $o(1/n^2)$ by choice of $A$. On the other hand, from Claim 15, for $A \log n \leq s \leq \sigma n$ we have

$$\mathbb{P}(|C_v \cap L| = s, |C_v \cap B| \geq |C_v \cap L|/2) \leq \exp(-cs) \leq \exp(-Ac \log n),$$

which we can again make $o(1/n^2)$ by choice of $A$. Hence, summing over $A \log n \leq s \leq \sigma n$, the probability that $v$ is in a small (size at most $\sigma n$, say) component of $H'$ containing at least $A \log n$ vertices of $L$ but not joined to $C_1$ in $H$ is $o(1/n)$, and the probability that such a vertex exists is $o(1)$. As noted above, this suffices to prove the theorem.

It remains to prove Claim 15. Recall that we have already conditioned on $H_R$; in particular, we have revealed all edges of $H'$ between vertices in $R$. It remains to reveal the edges of $H'$ between vertices in $L$, and between vertices in $L$ and vertices in $R \setminus C_1$. (By definition of $H'$, there are no edges of $H'$ between $L$ and $C_1$.) Let $v = v_1$ be a random vertex of $L$. We shall explore the component $C_v$ of $v$ in $H'$ in the following way: having "reached" vertices $v_1, \ldots, v_r \in C_v \cap L$ and "tested" $v_1, \ldots, v_{t-1}$, $t \leq r$, we next "test" vertex $v_t$. First, we add any neighbors (in the graph $H'$) of $v_t$ in $L$ not among the vertices reached so far to our list of reached vertices. Then we test edges between $v_t$ and $R \setminus C_1$, finding the set of components of $H_R$ that $v_t$ is joined to in $H'$. For each we find its unreached neighbors in $L$ and add them to our list.

The basic idea of the proof is simple, and similar to that of Lemma 14: roughly speaking, at each step we are much more likely to reach a good vertex of $L$ than a bad vertex. It will follow that the chance that $|C_v \cap L| = s$ but that the component $C_v$ contains at least $s/2$ bad vertices is exponentially small in $s$. The problem is that there is an exceptional case: when we reach an annoying component on the right, this may send many more edges in $H'$ to bad vertices than to good ones. But an annoying component of size $s'$ is very unlikely to have more than $s'$ bad neighbors, and the chance of reaching such a component will be exponentially small in $s'$, so the contribution from such annoying components is negligible. Turning this into a formal proof is now a matter of accounting.

Fix $s \leq \sigma n$, and recall that we must show that

$$\mathbb{P}(|C_v \cap L| = s, |C_v \cap B| \geq s/2) \leq \exp(-cs).$$

Let us define a quantity $f_t$ associated with the testing of vertex $v_t$: set

$$f_t = \exp(b - g/2 - 1/8)1_{E_t},$$

where $b$ and $g$ are the number of new good and bad vertices of $L$ that we reach when testing $v_t$, and $1_{E_t}$ is the indicator function of the event $E_t$ that after testing $v_t$ we have reached at most $s$ vertices in $C_v \cap L$. If $t > |C_v|$, so there is no vertex $v_t$ to test, set $f_t = 0$. The role of the indicator function $1_{E_t}$ is simply to stop our exploration process if we reach more than $s$ vertices in $C_v \cap L$, at which point there is nothing to prove.

Set $F_s = \prod_{t=1}^{s} f_t$. If $|C_v \cap L| = s$ and $|C_v \cap B| \geq s/2$, then

$$F_s = \exp(|(C_v \setminus \{v\}) \cap B| - |(C_v \setminus \{v\}) \cap (L \setminus B)|/2 - s/8)$$
$$\geq \exp(s/2 - 1 - s/4 - s/8) = \exp(s/8 - 1) \geq \exp(s/9),$$

assuming, as we may, that $s \geq 100$. Hence Claim 15 follows if we show that $\mathbb{E}(F_s) \leq 1$. In turn, it suffices to show that, conditional on the exploration so far, we have $\mathbb{E}(f_t) \leq 1$ for each $t$.

Let the small components of $H_R$ be $C_2, \ldots, C_m$. We shall test $v_t$ in several steps. In step 0 we check for edges from $v_t$ to unreached vertices in $L$. In step $i$, $2 \leq i \leq m$, if $C_i$ has not previously been reached by our exploration, we check for edges from $v_t$ to $C_i$, and, if we find such an edge, then test for edges from $C_i$ to unreached vertices in $L$. At every step we assume, as we may, that we have reached at most $s$ vertices in $L$; we shall suppress the corresponding indicator functions in the estimates below.

We may write $f_t$ as a product of a factor $f'_i$ (that also depends on $t$) for each step $i$; up to indicator functions corresponding to $1_{E_t}$, we may write

$$f'_0 = \exp(b_0 - g_0/2 - 1/8)$$

and
$$f'_i = \exp(b_i - g_i/2)$$
for $i = 2, \ldots, m$, where $b_i$ and $g_i$ are the number of new bad/good vertices reached in step $i$. (There is no step 1 as we do not test for edges to $C_1$.)

Since $|B| \leq 3\sigma n$ and, from (18), $\sigma < 1/(600\beta_{\max})$, for any $v \in V(G)$ we have $e_G(v, B) \leq n/200$. Hence, writing $B'$ for the set of vertices in $B$ not so far reached, conditional on everything so far, we have

$$\mathbb{E}(\exp(b_0)) = \prod_{w \in B'} (1 + (e-1)\beta_{v_t w}/n)$$

$$\leq \prod_{w \in B} \exp((e-1)\beta_{v_t w}/n)$$

(23)
$$= \exp((e-1)e_G(v_t, B)/n)$$

$$\leq \exp(1/100).$$

Using only $g_0 \geq 0$ it follows that $\mathbb{E}(f'_0) \leq \mathbb{E}(\exp(b_0)) \exp(-1/8) < \exp(-1/10)$.

Now let us condition not only on the results of testing $v_1, \ldots, v_{t-1}$, but also on steps $0, 2, \ldots, i-1$ of the testing of $v_t$, assuming as we may that we have reached at most $s \leq \sigma n$ vertices of $C_v \cap L$. Let $F_i$ be the event that we find an edge from $v_t$ to $C_i$.

Suppose that $C_i$ is not annoying. For any vertex $v \in R \setminus R^-$, the total weight of edges from $v$ to unreached vertices in $L$ is at least

(24) $\quad e_G(v, L) - \beta_{\max} s \geq 49\beta_{\max}\sigma n - \beta_{\max}\sigma n \geq 48\beta_{\max}\sigma n.$

On the other hand, for any vertex $v$,
$$e_G(v, B) \leq \beta_{\max}|B| \leq 3\beta_{\max}\sigma n.$$

As $|C_i \cap R^-| \leq |C_i|/2$, the total weight $w$ of edges of $G$ from $C_i$ to unreached vertices in $L \setminus B$ is at least $21\beta_{\max}\sigma n|C_i|$. A similar calculation to (23) but now using $1 - p(1 - e^{-1/2}) \leq \exp(-(1 - e^{-1/2})p) \leq \exp(-p/3)$ shows that $\mathbb{E}(\exp(-g_i/2) \mid F_i) \leq \exp(-7\beta_{\max}\sigma|C_i|)$. [Of course, we should update the set of reached vertices as we go, but (24) is valid at every step.] Arguing as for (23), since the total weight of edges of $G$ between $C_i$ and $B$ is at most $3\beta_{\max}\sigma n|C_i|$, we have $\mathbb{E}(\exp(b_i) \mid F_i, g_i) \leq \exp(6\beta_{\max}\sigma|C_i|)$, so $\mathbb{E}(f'_i \mid F_i) \leq 1$. Since $f'_i = 1$ whenever $F_i$ does not hold, it follows that the (conditional) expectation of $f'_i$ is at most 1.

Finally, suppose that $C_i$ is an annoying component of size $s'$. Using the bound $e_G(v, B) \leq \beta_{\max}|B| \leq 3\beta_{\max}\sigma n$ to bound $\mathbb{E}(\exp(b_i))$ as above, and using $g_i \geq 0$, we have $\mathbb{E}(f'_i \mid F_i) \leq \exp(6\beta_{\max}\sigma s')$. Since $\mathbb{P}(F_i) \leq \beta_{\max}|C_i|/n = \beta_{\max}s'/n$, it follows that

$$\mathbb{E}(f'_i) \leq 1 + \frac{\beta_{\max}s'}{n}(\exp(6\beta_{\max}\sigma s') - 1).$$

Since $\sigma \le 10^{-10}/\beta_{\max}^3$, for $s' \le 10^6 \beta_{\max}$, say, the bracket above is at most $1/(100\beta_{\max})$, so

$$\log(\mathbb{E}(f'_i)) \le \frac{s'}{100n}.$$

For $s' \ge 10^6 \beta_{\max}$ we use the very crude bound

$$\log(\mathbb{E}(f'_i)) \le \mathbb{E}(f'_i) - 1 \le \frac{\beta_{\max} s'}{n} \exp(6\beta_{\max}\sigma s').$$

Recall that $A_{s'}$ counts the number of vertices in annoying $s'$-vertex components, and that each such component contains $s'$ vertices. As all expectations are conditional on everything preceding them, we can multiply the expectations above together to conclude that

$$\log(\mathbb{E}(f_t)) \le -\frac{1}{10} + \sum_{s' \le 10^6 \beta_{\max}} \frac{A_{s'}}{100n}$$

$$+ \sum_{s' \ge 10^6 \beta_{\max}} \frac{\beta_{\max} A_{s'}}{n} \exp(6\beta_{\max}\sigma s').$$

By assumption, $A_{s'} \le \exp(-s'/200)n$ for every $s'$, while $\sum_i A_i \le n$. Hence

$$\log(\mathbb{E}(f_t)) \le -\frac{1}{10} + \frac{1}{100} + \sum_{s' \ge 10^6 \beta_{\max}} \beta_{\max} \exp(6\beta_{\max}\sigma s' - s'/200)$$

$$\le -0.09 + \sum_{s' \ge 10^6 \beta_{\max}} \beta_{\max} \exp(-s'/400) < 0.$$

In other words, $\mathbb{E}(f_t) < 1$. Recalling that the argument above, and hence the final estimate, hold conditional on all previous steps in the exploration, it follows that

$$\mathbb{E}(F_s) = \mathbb{E}\left(\prod_{i=1}^s f_t\right) < 1.$$

As noted earlier, this implies Claim 15 and hence Theorem 9.  □

The proof of Theorem 9 presented above is rather involved, and the reader may well wonder whether it can be simplified. While writing this paper, we found various proofs of Theorem 9 that were indeed much simpler; unfortunately they were also incorrect. Of course, there may well be a simple proof that we have missed. If so, it would be interesting to find one.

**Acknowledgment.** The authors would like to thank Svante Janson for a very careful reading of the first version of this paper, leading to several minor corrections.

B. BOLLOBÁS  
DEPARTMENT OF PURE MATHEMATICS  
  AND MATHEMATICAL STATISTICS  
UNIVERSITY OF CAMBRIDGE  
WILBERFORCE ROAD  
CAMBRIDGE CB3 0WB  
UNITED KINGDOM  
AND  
DEPARTMENT OF MATHEMATICAL SCIENCES  
UNIVERSITY OF MEMPHIS  
MEMPHIS, TENNESSEE 38152  
USA  
E-MAIL: b.bollobas@dpmms.cam.ac.uk  

C. BORGS  
J. CHAYES  
MICROSOFT RESEARCH NEW ENGLAND  
1 MEMORIAL DRIVE  
CAMBRIDGE, MASSACHUSETTS 01242  
USA  
E-MAIL: borgs@microsoft.com  
        jchayes@microsoft.com  

O. RIORDAN  
MATHEMATICAL INSTITUTE  
UNIVERSITY OF OXFORD  
24–29 ST GILES'  
OXFORD OX1 3LB  
UNITED KINGDOM  
E-MAIL: riordan@maths.ox.ac.uk